\documentclass[12pt,reqno]{amsart}
\usepackage[toc,page]{appendix}

\makeatletter
\def\@tocline#1#2#3#4#5#6#7{\relax
  \ifnum #1>\c@tocdepth 
  \else
    \par \addpenalty\@secpenalty\addvspace{#2}%
    \begingroup \hyphenpenalty\@M
    \@ifempty{#4}{%
      \@tempdima\csname r@tocindent\number#1\endcsname\relax
    }{%
      \@tempdima#4\relax
    }%
    \parindent\z@ \leftskip#3\relax \advance\leftskip\@tempdima\relax
    \rightskip\@pnumwidth plus4em \parfillskip-\@pnumwidth
    #5\leavevmode\hskip-\@tempdima
      \ifcase #1
       \or\or \hskip 1em \or \hskip 2em \else \hskip 3em \fi%
      #6\nobreak\relax
      \dotfill
      \hbox to\@pnumwidth{\@tocpagenum{#7}}
    \par
    \nobreak
    \endgroup
  \fi}
\makeatother

\usepackage{physics}

\usepackage{amssymb,amsthm,amsmath,amsxtra,dsfont,appendix,bbm,stix,xparse}

\usepackage{hyperref} 
\usepackage{comment}
\numberwithin{equation}{section}

\usepackage{graphicx} 

\theoremstyle{plain}
\newtheorem{theorem}{Theorem}[section]
\newtheorem{cor}[theorem]{Corollary}
\newtheorem{lemma}[theorem]{Lemma}

\newtheorem{prop}[theorem]{Proposition}
\newtheorem{assumption}[theorem]{Assumption}

\theoremstyle{definition}
\newtheorem{definition}[theorem]{Definition}

\usepackage{colonequals} 

\usepackage{color}

\newcommand{\LR}[1]{\left(#1\right)}

\title{Dimensional reduction for anyons in the average-field approximation}

\author[Q. Yang]{Qiyun Yang}
\address{Ecole Normale Sup\'erieure de Lyon, UMPA (UMR 5669)}
\email{qiyun.yang@ens-lyon.fr}

\date{May 2026}

\begin{document}


\maketitle
%

\begin{abstract}
    We study abelian anyons at the mean-field/almost-bosonic level, whose dynamics are governed by the Chern-Simons-Schr\"odinger system. We consider the dimensional reduction of this 2D model by introducing an anisotropic trapping potential, and derive an effective 1D model after tracing out the tight confinement direction. The resulting effective dynamics in the loose confinement direction is captured by a quintic defocusing nonlinear Schr\"odinger equation. We rigorously establish this dimensional reduction process in the sense of ground state energies and time-dependent solutions, under a uniform $H^2$ well-posedness assumption.
\end{abstract}


\tableofcontents


\section{Introduction}

\subsection{Motivation}

Anyons, quasi-particles with fractional statistics interpolating between bosons and fermions, are important objects in two-dimensional quantum physics. 
Their theoretical description often relies on the Chern–Simons–Schr\"odinger (CSS) equation, which provides an effective model for the emergence of exotic exchange phases and plays an important role in the fractional quantum Hall effect (see e.g. \cite[Section 5.16]{Jain-07}\cite{Khare-05,LopFra-91,Wilczek-90,Zhang-92,ZhaHanKiv-89}).
Since CSS theory couples a matter field to a gauge field, a promising research direction is to realize aspects of 
anyon physics in the setting of cold atoms interacting with artificial gauge fields (see e.g. \cite{ChiEtalCel-22,DalGerJuzOhb-11,EdmEtalOhb-13,FroEtalTar-22,IacCabTarCel-25,LunRou-16,ValOhb-24,ValWesOhb-20,ZhaSreGemJai-14,ZhaSreJai-15}).

Understanding dimensional reduction is particularly relevant for connecting 2D anyon models to experimentally feasible quasi 1D systems. 
A strong confinement in one spatial direction provides a natural mechanism for dimensional reduction, which can be implemented in cold-atom experiments.
Over the past decades, analogous confinement-induced dimensional reductions have been extensively investigated in various particle systems, particularly in Bose–Einstein condensates (see e.g. \cite[Chapter 8\&9]{LieSeiSolYng-05}\cite{AftBla-08,BaoTreMeh-15,BaoTreMeh-17,Ben-Cas-Meh,BenMehSchWei-05,MehSpa-16} and references therein).

In previous work, we analyzed the dimensional reduction of the 2D many-body abelian anyon model, where the resulting 1D limit in the loosely confined direction is described by the Tonks–Girardeau (TG) gas \cite{RouYan-23a,RouYan-23b}. 
Since the mean-field (almost-bosonic) limit of the 2D many-body anyons is governed by the CSS equation (see e.g. \cite{AtaEllGetGirLunNgu-26,AtaGirLun-25,ChiSen-92,CorDubLunRou-19,Girardot-19,Girardot-25b,GirLee-24,LunRou-15,Visconti-25} and references therein), it is natural to consider the problem of commuting these two limiting processes. 
For further physics discussions, we refer readers to \cite{RouYan-26}.

\subsection{Model}

We consider the 2D model of abelian anyons in the average-field approximation described by the energy functional (Hamiltonian)
\begin{equation} 
    \mathcal{E}^{2\mathrm{D}}_{\varepsilon} (\psi) = \int_{\mathbb{R}^2 } \left| \left( - \mathrm{i} \nabla_{} +  \beta \mathbf{A}[|\psi|^2]\right) \psi  \right|^2  + \int_{\mathbb{R}^2} V_{\varepsilon} |\psi|^2,  \label{2Denergy}
\end{equation}
where $\beta \in \mathbb{R}$ is an effective coupling constant, $V_{\varepsilon}$ is the anisotropic trapping potential 
\begin{equation}\label{def-v}
    V_{\varepsilon}(x,y) = |x|^2 + \varepsilon^{-2} |y|^2 \quad \text{for} \quad 0 < \varepsilon \ll  1
\end{equation}
and $\mathbf{A[\cdot]}$ is the density-dependent vector potential in Coulomb gauge:
\begin{equation}\nonumber
    \mathbf{A}[\rho](\mathbf{x})  =  \int_{\mathbb{R}^2} \frac{(\mathbf{x}-\mathbf{x}')^{\perp}}{|\mathbf{x} - \mathbf{x}'|^2} \rho(\mathbf{x}') \mathrm{d} \mathbf{x}' = \left( (\nabla^{\perp}_{\mathbf{x}} \omega_0) * \rho \right) (\mathbf{x}) 
\end{equation}
for
\begin{equation}\nonumber
    \mathbf{x} = (x,y), \quad  \mathbf{x}^{\perp} = (-y ,x), \quad \nabla^{\perp}=(-\partial_y, \partial_x), \quad  \omega_0(\mathbf{x}) = \log |\mathbf{x}|.
\end{equation}
The corresponding 2D dynamics for \eqref{2Denergy} is governed by the Chern-Simons-Schr\"odinger (CSS) equation \footnote{One can refer to \cite[Lemma A.2]{CorLunRou-17} for detailed calculations of the variational equation.} \footnote{Equation \eqref{2DPDE} can be derived from the Lagrangian for a CSS system where both $\mathbf{A}$ and $\psi$ are variables. See e.g. \cite{BerBouSau-95,LiuSmiTat-14,RouYan-26}.}
\begin{equation} \label{2DPDE}
    \mathrm{i}  \partial_t \psi  = \delta_{\overline{\psi}} \mathcal{E}^{2\mathrm{D}}_{\varepsilon} (\psi)  = \left[ \left( -\mathrm{i}\nabla_{\mathbf{x}}+\beta \mathbf{A}[|\psi|^2] \right)^2  - 2\beta (\nabla_{\mathbf{x}}^{\perp} \omega_0) * \mathbf{J}_{\beta \mathbf{A} [|\psi|^2]}(\psi) \right] \psi + V_{\varepsilon} \psi, 
\end{equation}
where 
\begin{align}\label{defJ}
    \mathbf{J}_{\mathbf{A}}(\psi) & = \frac{1}{2} \left[ \overline{\psi} \left( -\mathrm{i} \nabla_{\mathbf{x}} + \mathbf{A} \right) \psi  +  {\psi} \overline{\left( -\mathrm{i} \nabla_{\mathbf{x}} + \mathbf{A} \right) \psi } \right] 
\end{align}
is the current.
The model introduced above indeed gives an effective mean-field description of many-body abelian anyons around the bosonic end (almost bosonic anyons). 
The effective coupling constant $\beta$ still depends on the exchange parameter $\alpha$ for many-body anyons following the relation 
\begin{equation}\nonumber
    \beta \sim \alpha N,
\end{equation}
where $N$ denotes the number of particles.
In this mean-field limit, $\beta$ is kept fixed while $N$ tends to infinity, implying that $\alpha$ must approach 0, corresponding to the almost-bosonic anyons regime.

Since particles are more likely to stay at a lower energy, they will gather around the one-dimensional line $y=0$ as the parameter $\varepsilon \to 0$ due to the anisotropic trapping potential $V_{\varepsilon}$ as in \eqref{def-v}. 
In the following sections, we prove that in both static (i.e. for ground state) and dynamical (with a uniform $H^2$ well-posedness assumption) aspects, the effective 1D limit model on the line $y=0$ as $\varepsilon$ goes to 0 is described by the energy functional (Hamiltonian)
\begin{equation}\label{1Denergy}
    \mathcal{E}^{1\mathrm{D}}_{} (\varphi) = \int_{\mathbb{R}} \left( |\partial_x \varphi (x)|^2  + \frac{1}{3} \pi^2 \beta^2 |\varphi(x)|^6   +  |x|^2 |\varphi(x)|^2 \right) \mathrm{d} x,
\end{equation}
whose corresponding 1D dynamics is the quintic nonlinear Schr\"odinger (quintic NLS) equation
\begin{equation}\label{1dPDE}
    \mathrm{i} \partial_t \varphi  = \delta_{\overline{\varphi}} \mathcal{E}^{1\mathrm{D}}_{} (\varphi) = -\partial_x^2 \varphi + \pi^2 \beta^2 |\varphi|^4 \varphi + |x|^2 \varphi.
\end{equation}
In the $y$-direction, the harmonic part dominates the energy. We let
\begin{equation}\label{def-u}
    u_{\varepsilon}(y) = \left( \sqrt{\pi \varepsilon} \right)^{-\frac{1}{2}} e^{-\frac{y^2}{2\varepsilon}} 
\end{equation}
be the ground state of the harmonic oscillator 
$$-\partial^2_y + \varepsilon^{-2} |y|^2,$$ 
and its corresponding energy is 
\begin{equation}\label{OHenergy}
    e_{\varepsilon} = \frac{1}{\varepsilon}.
\end{equation}

The correct ansatz for the ground state ($L^2$ normalized minimizer) of $\mathcal{E}^{2\mathrm{D}}_{\varepsilon}$ almost looks like
\begin{equation}\label{ansatz}
    \varphi(x) u_{\varepsilon}(y) e^{-\mathrm{i} \beta S[|\varphi|^2u_{\varepsilon}^2](x,y)}
\end{equation}
for 
\begin{equation}   \label{defS}
    S[\rho](\mathbf{x}) = \int_{\mathbb{R}^2} S(\mathbf{x}-\mathbf{y}) \rho(\mathbf{y}) \mathrm{d} \mathbf{y} = (S*\rho) (\mathbf{x}),  \quad \quad S(x,y) = \arctan \frac{y}{x}. 
\end{equation}
The reason why we introduce the phase factor $S[\rho]$ is that its gradient is related to $\mathbf{A}[\rho]$:
\begin{equation}\nonumber
    \nabla (S[\rho]) = (\nabla S)*\rho = \mathbf{A}[\rho]  - \mathbf{T}  [\rho] ,
\end{equation}
where 
\begin{equation}\label{def-T}
    \mathbf{T} [\rho] (x,y) = \LR{\begin{matrix}
        - \pi \int_{\mathbb{R}}\operatorname{sgn}(y-y') \rho (x,y') \differential y'\\
        0
        
    \end{matrix}} = (\mathbf{T}_0 * \rho) (x,y)
\end{equation} 
with
\begin{equation}\nonumber
    \mathbf{T}_0(x,y) = 
    \left( \begin{matrix}
        -\pi \operatorname{sgn}(y)\delta_{x=0} \\
        0
    \end{matrix}\right) 
    \quad \text{ for } \quad 
    \operatorname{sgn}(y) = \begin{cases}
        1 & \text{ when } y>0\\
        0 & \text{ when } y=0 \\
        -1 & \text{ when } y<0
    \end{cases}
\end{equation}
and $\delta_{x=0}$ is the Dirac delta function (distribution). 
Notice that $\mathbf{T}[\rho]$ has 0 component in the $y$ direction, which makes it more convenient than $\mathbf{A}[\rho]$ for dimensional reduction analysis.
But the cost is that $\mathbf{T}[\rho] $ has less regularity than $\mathbf{A}[\rho]$ given the same $\rho$.
Thanks to the phase factor $S[\rho]$ and a change of gauge, we can simplify the Coulomb gauge $\mathbf{A}[\rho]$ to the lineal gauge $\mathbf{T}[\rho]$.
The calculations in Section \ref{secUpper} show that the 2D energy of the ansatz \eqref{ansatz} precisely equals the sum of $e_\varepsilon$, the ground state energy in the $y$ direction, and the 1D quintic NLS energy of $\varphi$.

\subsection{Results}

We rigorously prove the dimensional reduction process mentioned above at the ground state level (Theorem \ref{GSEresult} and Theorem \ref{gsfunction}) and at the time evolution level under a uniform $H^2$ well-posedness assumption (Theorem \ref{dynamicresult}).

\begin{definition}[Ground state energies] \label{GSE}
    We denote $E^{2\mathrm{D}}_{\varepsilon}$ by the ground state energy of the 2D anyonic model described by \eqref{2Denergy}, which is defined as
    \begin{equation}\nonumber
        E^{2\mathrm{D}}_{\varepsilon} = \inf \left\{ \mathcal{E}^{2\mathrm{D}}_{\varepsilon} (\psi) : \norm{\psi}_{L^2(\mathbb{R}^{2})} = 1 \right\}.
    \end{equation}
    We denote $E^{1\mathrm{D}}_{}$ by the ground state energy of the 1D quintic NLS model described by \eqref{1Denergy}, which is defined as
    \begin{equation}\nonumber
        E^{1\mathrm{D}}_{} = \inf \left\{ \mathcal{E}^{1\mathrm{D}}_{} (\varphi) : \norm{\varphi}_{L^2(\mathbb{R}^{})} = 1 \right\}.
    \end{equation}
\end{definition}

\begin{theorem}[\bf{Ground state energies}] \label{GSEresult}
Let $E^{2\mathrm{D}}_{\varepsilon} $ and $E^{1\mathrm{D}}_{} $ be as in Definition \ref{GSE}, and let $e_{\varepsilon}$ be as in \eqref{OHenergy}. Then we have 
\begin{equation}
    \lim_{\varepsilon \to 0} \left( {E}^{2\mathrm{D}}_{\varepsilon} - e_{\varepsilon} \right) = {E}^{1\mathrm{D}}_{} .  \nonumber
\end{equation}
\end{theorem}
\begin{theorem}[\bf Ground states / $L^2$ normalized minimizers]\label{gsfunction}
    Let $\Psi$ be a ground state for \eqref{2Denergy} and let $u_{\varepsilon}$ be as in \eqref{def-u}. Then there exists\footnote{$\varphi_0$ is unique up to a phase. Different vanishing sequence $(\varepsilon_k)_{k\in \mathbb{N}}$ might result in $\varphi_0$ with different phase.} a ground state $\varphi_0$ for \eqref{1Denergy} and a vanishing sequence $(\varepsilon_k)_{k\in \mathbb{N}}$ such that
    \begin{equation}\nonumber
        \lim_{k\to \infty} \norm{\Psi - \varphi_0 u_{\varepsilon_k}}_{L^2{(\mathbb{R}^2 )} } = 0.
    \end{equation}
\end{theorem}
Theorem \ref{GSEresult} indicates that, in this dimensional reduction process, the ground state energy of 2D anyons is asymptotically equal to the sum of the ground state energies in two directions, the tight one in the $y$ direction and the loose one in the $x$ direction, where the ground state energy in the $x$ direction is given by the quintic NLS energy. 
And the 2D ground state also decoupled into these two directions as stated in Theorem \ref{gsfunction}.
The phase $e^{\mathrm{i} \beta S}$ in \eqref{ansatz} in not visible at the $L^2$ level, but it would appear in an $H^1$ bound, as in the decoupling relation for many anyons proved in \cite{RouYan-23b}.

\begin{definition}\label{defSigma}
    For $n=1,2$ and $s=1,2$, we define \footnote{Readers can find in Appendix \ref{equivnorm} the norm on this space.}
\begin{equation}\nonumber
    \Sigma^s (\mathbb{R}^n) \colonequals \left\{ \Psi \in L^2(\mathbb{R}^n) : (-\Delta + |\mathbf{x}|^2)^{\frac{s}{2}}\Psi \in L^2(\mathbb{R}^n) \right\}.
\end{equation}
\end{definition} 

Let $\varphi_0 \in \Sigma^2 (\mathbb{R})$ be normalized in $L^2(\mathbb{R})$ and $\mathcal{E}^{1\mathrm{D}}(\varphi_0) < \infty $. 
Consider the initial value problem (IVP) for the 2D anyonic dynamics \eqref{2DPDE}:
\begin{equation}\label{IVP-A}
\begin{cases}
    \mathrm{i}  \partial_t \psi  =  \big[ \left( -\mathrm{i}\nabla_{\mathbf{x}} + \beta \mathbf{A}[|\psi|^2] \right)^2  - 2\beta (\nabla_{\mathbf{x}}^{\perp} \omega_0) * \mathbf{J}_{\beta \mathbf{A} [|\psi|^2]}(\psi) \big] \psi + V_{\varepsilon} \psi\\
    \psi(0,x,y) = \varphi_0(x) u_{\varepsilon}(y) e^{- \mathrm{i }\beta S[|\varphi_0|^2 u_{\varepsilon}^2](x,y)}
\end{cases}.
\end{equation} 
The choice of this initial datum is motivated by the ansatz \eqref{ansatz}.
For simplicity, we perform a change of gauge and a rescaling on the original IVP \eqref{IVP-A} before taking $\varepsilon$ to 0. 
In order to change the Coulomb gauge $\mathbf{A}[\rho]$ to the lineal gauge $\mathbf{T}[\rho]$, we rewrite the solution $\psi$ to the original IVP \eqref{IVP-A} as 
\begin{equation} \nonumber
    \psi (t, x, y) = \tilde{\psi} (t,x,y) e^{-\mathrm{i} \beta S[|\tilde{\psi}|^2](t,x,y)} \quad \text{for} \quad \tilde{\psi} (t,x,y) = \psi (t, x, y)  e^{\mathrm{i} \beta S[|\psi|^2](t,x,y)},
\end{equation}
then the two-dimensional IVP \eqref{IVP-A} becomes 
\begin{equation} \label{gaugedawayPDE}
\begin{cases}
    \mathrm{i}  \partial_t \tilde{\psi}  = \big[ \left( -\mathrm{i}\nabla_{\mathbf{x}}+\beta \mathbf{T}[|\tilde{\psi}|^2] \right)^2  - 2\beta \mathbf{T}_0 * \mathbf{J}_{\beta \mathbf{T} [|\tilde{\psi}|^2]}(\tilde{\psi}) \big] \tilde{\psi} + V_{\varepsilon} \tilde{\psi}\\
    \tilde{\psi} (0,x,y) = \varphi_0(x) u_{\varepsilon}(y) 
\end{cases}.
\end{equation}
To balance the energies in two directions, we rescale the solution $\tilde{\psi}$ to the gauged IVP \eqref{gaugedawayPDE} as
\begin{equation}\nonumber
    \phi_{\varepsilon} (t,x,y) = {\varepsilon}^{\frac{1}{4}} \tilde{\psi} \left(t,x, \sqrt{\varepsilon}y\right) 
\end{equation}
to obtain 
\begin{equation}\nonumber
\begin{cases}
    \mathrm{i}  \partial_t \phi_{\varepsilon}   =  \big[ \left( -\mathrm{i}\nabla_{\mathbf{x},\varepsilon}+\beta \mathbf{T}[|\phi_{\varepsilon}|^2] \right)^2  - 2\beta \mathbf{T}_0 * \mathbf{J}_{\beta \mathbf{T} [|\phi_{\varepsilon}|^2]}(\phi_{\varepsilon}) \big] \phi_{\varepsilon} + V_{\sqrt{\varepsilon}} \phi_{\varepsilon}\\
    \phi_{\varepsilon}(0,x,y) =  \varphi_0(x) u_{1}(y) 
\end{cases},
\end{equation}
which can be rewritten as
\begin{equation} \label{rescaledPDE}
\begin{cases}
    \mathrm{i}  \partial_t \phi_{\varepsilon}   = \frac{1}{\varepsilon} H_y \phi_{\varepsilon} + H_x \phi_{\varepsilon} + f[\phi_{\varepsilon}]\\
    \phi_{\varepsilon}(0,x,y) =  \varphi_0(x) u_{1}(y) 
\end{cases}
\end{equation}
with 
\begin{equation}\label{H_xH_y}
    H_x = - \partial_x^2 +|x|^2, \quad   H_y = - \partial_y^2 +|y|^2
\end{equation}
and
\begin{multline} \label{f[]}
    f[\phi_{\varepsilon}] = \beta^2 \left( (\mathbf{T}_0)_x *|\phi_{\varepsilon}|^2 \right)^2 \phi_{\varepsilon} - \mathrm{i} \beta \left( (\mathbf{T}_0)_x *|\phi_{\varepsilon}|^2  \right) \partial_x \phi_{\varepsilon}   - \mathrm{i} \beta \partial_x \left( \left( (\mathbf{T}_0)_x *|\phi_{\varepsilon}|^2 \right) \phi_{\varepsilon}   \right) \\ - 2 \beta (\mathbf{T}_0)_x * \left( (\mathbf{J}_{\mathbf{0}}(\phi_{\varepsilon}) )_x +\beta \left( (\mathbf{T}_0)_x *|\phi_{\varepsilon}|^2  \right) |\phi_{\varepsilon}|^2 \right) \phi_{\varepsilon},
\end{multline}
where  
$$ \nabla_{\mathbf{x},\varepsilon} = \LR{\partial_x , {\varepsilon}^{-\frac{1}{2}} \partial_y},$$
$(\mathbf{T}_0)_x$ is defined as the $x$ component of $\mathbf{T}_0$, i.e.
\begin{equation}\label{T0x}
    (\mathbf{T}_0)_x (x,y) = - \pi \operatorname{sgn} (y) \delta_{x=0},
\end{equation}
and similarly for $(\mathbf{J}_{\mathbf{0}}(\phi_{\varepsilon}) )_x $.

\begin{assumption}\label{assumption}
There exists a constant $T_0 >0$ depending on the initial datum $\varphi_0$, and a constant $\varepsilon_{{0}} >0 $ depending on $T_0$ such that the initial value problem \eqref{rescaledPDE} has a unique solution $\phi_{\varepsilon} \in \mathcal{C} ( [0,T_0], \Sigma^2(\mathbb{R}^2))$, which is uniformly bounded with respect to $\varepsilon \in ( 0, \varepsilon_0]$ in $\mathcal{C} ( [0,T_0], \Sigma^2(\mathbb{R}^2))$.
\end{assumption}
To achieve the proof for the dynamics, we need the above assumption about uniform $H^2$ well-posedness on the IVP \eqref{rescaledPDE}.
The uniform boundedness assumption plays an important role in the essential estimate stated in Proposition \ref{estimateh}.
In view of the results and methods in \cite{BerBouSau-95}, we have the uniform $H^2$ well-posedness with $T_0 \sim \varepsilon$ for the rescaled IVP
\begin{equation}\label{IVP-Aresacled}
\begin{cases}
    \mathrm{i}  \partial_t \psi_{\varepsilon}  =  \big[ \left( -\mathrm{i}\nabla_{\mathbf{x},\varepsilon} + \beta \mathbf{A}_{\varepsilon}[|\psi_{\varepsilon}|^2] \right)^2  - 2\beta (\nabla_{\mathbf{x}}^{\perp} \omega_0)_{\varepsilon} * \mathbf{J}_{\beta \mathbf{A}_{\varepsilon} [|\psi_{\varepsilon}|^2],\varepsilon}(\psi_{\varepsilon}) \big] \psi_{\varepsilon} + V_{\sqrt{\varepsilon}} \psi_{\varepsilon}\\
    \psi_{\varepsilon}(0,x,y) = \varphi_0(x) u_{1}(y) e^{- \mathrm{i }\beta S_{\varepsilon}[|\varphi_0|^2 u_{1}^2](x,y)}
\end{cases},
\end{equation} 
where 
$$ \psi_{\varepsilon}(t,x,y) = \psi \LR{t,x,\sqrt{\varepsilon}y} $$
$$ \mathbf{A}_{\varepsilon}[\rho] = (\nabla_{\mathbf{x}}^{\perp} \omega_0)_{\varepsilon} * \rho, \quad (\nabla_{\mathbf{x}}^{\perp} \omega_0)_{\varepsilon} (x,y) = (\nabla_{\mathbf{x}}^{\perp} \omega_0 )\LR{x, \sqrt{\varepsilon} y} , $$
$$ \mathbf{J}_{\mathbf{A},\varepsilon}(\psi) = \frac{1}{2} \left[ \overline{\psi} \left( -\mathrm{i} \nabla_{\mathbf{x},\varepsilon} + \mathbf{A} \right) \psi  +  {\psi} \overline{\left( -\mathrm{i} \nabla_{\mathbf{x},\varepsilon} + \mathbf{A} \right) \psi } \right], $$
\begin{equation} \label{def-Sepsilon}
    S_{\varepsilon} [\rho] = S_{\varepsilon} * \rho, \quad S_{\varepsilon}(x,y) = S\LR{x, \sqrt{\varepsilon} y } .
\end{equation}
The uniform $H^2$ well-posedness we stated in Assumption \ref{assumption} is equivalent to the corresponding uniform $H^2$ well-posedness on the rescaled IVP \eqref{IVP-Aresacled}.
The difficulty of improving the evolution time $T_0$ arises from the rescaled gauge $\mathbf{A}_{\varepsilon}[\rho]$, which depends on $\varepsilon$.
On the other hand, the kernel $\mathbf{T}_0$ of the lineal gauge $\mathbf{T}[\rho]$ doesn't change after the rescaling on $y$ direction, hence the gauge in the IVP \eqref{rescaledPDE} remain the same lineal gauge $\mathbf{T}[\rho]$ independent of $\varepsilon$.
But the fact that $\mathbf{T}[\rho]$ has less regularity than $\mathbf{A}_{\varepsilon}[\rho]$ troubled us a lot when trying to prove the $H^2$ well-posedness by Kato's method \cite[Theorem 6]{Kato-75}: for example, estimating the commutator $[\partial_x^2, (\mathbf{T}_0)_x * (\mathbf{J}_{\mathbf{0}}(\phi_{\varepsilon}) )_x ]$ requires integrability of $\partial_x^3 \phi_{\varepsilon}$, which is unclear.
For more related work on well-posedness of the CSS equation, see, e.g., \cite{Ataei-24,Ataei-25,AtaLunNgu-24,GirLee-24,Lim-18,LiuSmiTat-14,OhPus-15}.
Ataei proved in \cite[Theorem 1.1]{Ataei-24} a sharp condition on $H^1$ global well-posedness of the CSS equation, and Lim proved in \cite[Theorem 1.1]{Lim-18} the criterion that the global $H^s$ solution with $s\ge 1$ exists if and only if the global $H^1$ solution exists.
With these two results, the CSS equation with zero $\phi^4$ interaction term and zero external potential has a unique global $H^2$ solution.
We believe that the external harmonic potential will not change the existence of the solution, but the uniform boundedness of the solutions with respect to $\varepsilon$ is a question.
However, our dynamic result below (Theorem \ref{dynamicresult}) is also true for evolution time $T_0 \sim \varepsilon$. 
Apart from that, Assumption \ref{assumption} is true if the coupling constant $\beta \sim \varepsilon$, but this case results in a free harmonic oscillator instead of the quintic NLS in the limit, which is quite trivial.

\begin{theorem}[\bf{Dynamics}]\label{dynamicresult}
Let $\psi$ and $\phi_{\varepsilon}$ be the solutions to the initial value problems \eqref{IVP-A} and \eqref{rescaledPDE} respectively, and let $\varphi$ be the solution to
\begin{equation} \label{1Divp}
    \begin{cases}
        \mathrm{i} \partial_t \varphi = -\partial_x^2 \varphi + |x|^2 \varphi + \pi^2 \beta^2 |\varphi|^4 \varphi \\
        \varphi(0,x) = \varphi_0 (x)
    \end{cases}.
\end{equation}
Under Assumption \ref{assumption} above, there exists a constant $C_{T_0}$ independent of $\varepsilon$ such that
    \begin{equation}\nonumber
        \sup_{t\in (0,T_{0})} \norm{\phi_{\varepsilon}(t,\cdot) - e^{-\mathrm{i}t \varepsilon^{-1} } \varphi(t,\cdot)u_{1} }_{L^2} \le C_{T_0} \varepsilon^{\frac{1}{4}},
    \end{equation}
    \begin{equation}\nonumber
        \sup_{t\in (0,T_{0})} \norm{\psi(t,\cdot) - e^{-\mathrm{i}t \varepsilon^{-1} } \varphi(t,\cdot)u_{\varepsilon} }_{L^2} \xrightarrow{\varepsilon \to 0} 0.
    \end{equation}
\end{theorem}

Theorem \ref{dynamicresult} shows that, in this dimensional reduction process, the dynamics of the 2D anyons also decomposes into the two directions, given a proper initial datum.

\begin{figure}[htp]
    \centering
    \includegraphics[width=0.8\linewidth]{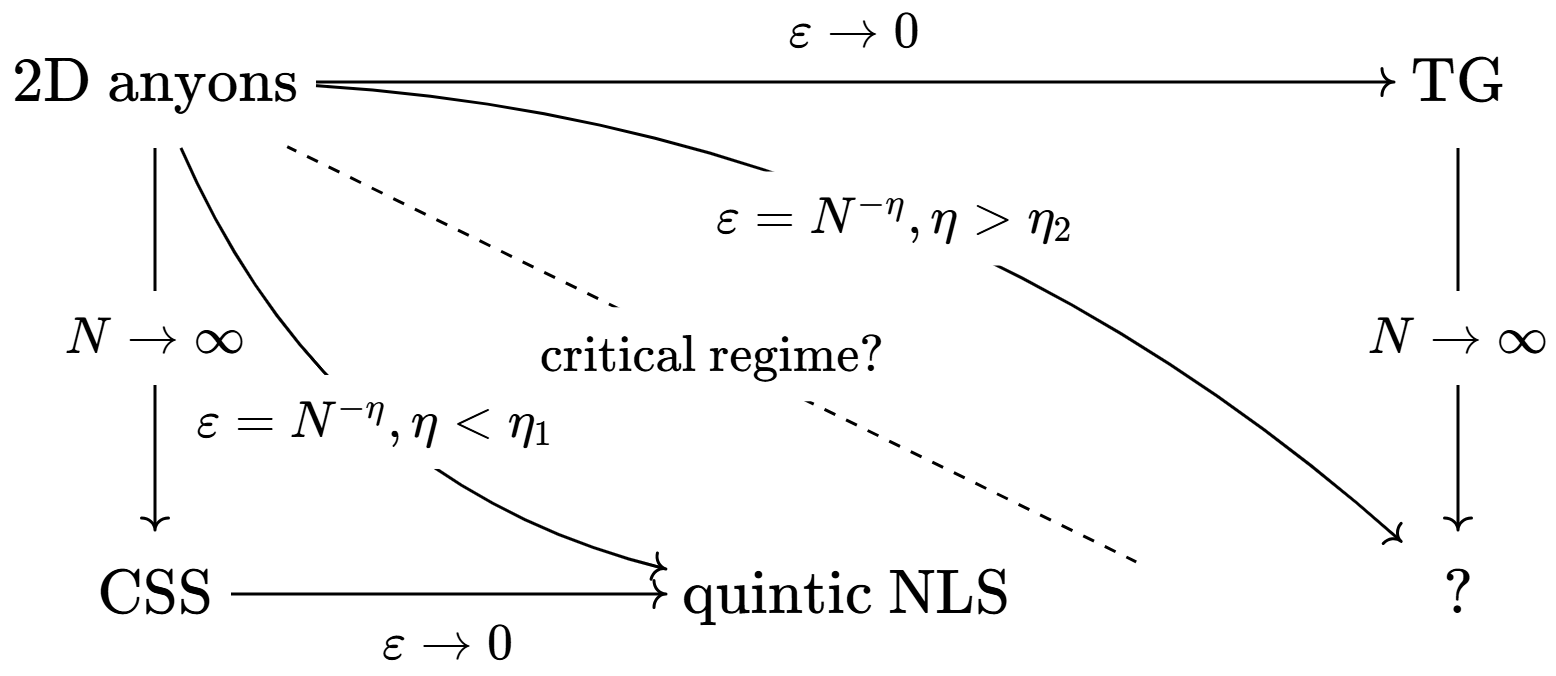}
    \caption{Exchange diagram for different limit processes on anyons.}
    \label{exchangediam}
\end{figure}

The quintic NLS equation \eqref{1dPDE} obtained as the 1D limit model still depends on the exchange parameter $\alpha$ through the coupling constant $\beta$, whereas the TG gas, which arises as the 1D limit of the many-body anyons model \cite{RouYan-23a,RouYan-23b}, does not.
This indicates that the mean-field limit does not commute with the dimensional reduction limit as illustrated in Figure \ref{exchangediam}. 
There may exist a critical regime between these two limiting processes, which remains an open question.

The proof of Theorem \ref{GSEresult} is divided into an energy upper bound part and an energy lower bound part in Section \ref{GSEproof}, and Theorem \ref{gsfunction} results from the proof of the energy lower bound and Section \ref{SecL2cvg}. 
A heuristic calculation is given in Section \ref{trialPDE} to make Theorem \ref{dynamicresult} plausible, and in Section \ref{dynamicsproof}, the full proof for Theorem \ref{dynamicresult} is given.\\

\noindent\textbf{Acknowledgments.} 
Thanks to Nicolas Rougerie for discussions, proofreading, and the joint work on \cite{RouYan-26}.
Thanks to Th\'eotime Girardot for discussions about the CSS system. 
Thanks to Chengzhi Ye for discussions on various physics topics.

\section{Dimensional reduction for ground state energies and ground states} \label{GSEproof}

In this section, we investigate the dimensional reduction process for ground state energies and ground states, from the CSS equation to the quintic NLS equation. 
Lemma \ref{upperbound} in Section \ref{secUpper} and Lemma \ref{lowerbound} in Section \ref{secLower} result in the relation between ground state energies stated in Theorem \ref{GSEresult}.
The proof for Theorem \ref{gsfunction} ends in Section \ref{SecL2cvg}.

\subsection{Energy upper bound}\label{secUpper}

In this section, we prove the following:

\begin{lemma}[Energy upper bound] \label{upperbound}
Let $E^{2\mathrm{D}}_{\varepsilon} $ and $E^{1\mathrm{D}}_{} $ be as in Definition \ref{GSE}, and let $e_{\varepsilon}$ be as in \eqref{OHenergy}. Then we have
\begin{equation}\nonumber
     {E}^{2\mathrm{D}}_{\varepsilon} \le  e_{\varepsilon}  +  {E}^{1\mathrm{D}}_{}.
\end{equation}
\end{lemma}

\begin{proof}
    Consider a trial state
\begin{equation} \nonumber
    \psi(x,y) = \varphi(x) u_{\varepsilon}(y) e^{-\mathrm{i}\beta S[|\varphi |^2 u_{\varepsilon}^2](x,y)},
\end{equation}
where $u_{\varepsilon}$ is as in \eqref{def-u}, $S[\cdot]$ is as in \eqref{defS}, and $\varphi$ to be determined is normalized in $L^2$.
Performing the change of gauge, we have
\begin{equation}\label{get-rid-of-phase}
     \left( - \mathrm{i} \nabla_{} +  \beta \mathbf{A}[|\psi|^2]\right) \psi  =   e^{-\mathrm{i}\beta S[|\varphi u_{\varepsilon}|^2]} \left( - \mathrm{i} \nabla_{} +  \beta \mathbf{T} [|\varphi u_{\varepsilon}|^2] \right) (\varphi u_{\varepsilon })  ,
\end{equation}
where $\mathbf{T}[\cdot]$ is defined as in \eqref{def-T}.
And the 2D energy \eqref{2Denergy} becomes
\begin{align}
    \mathcal{E}^{2\mathrm{D}}_{\varepsilon} (\psi) & = \int_{\mathbb{R}^2 } \left| \left( - \mathrm{i} \nabla_{} +  \beta \mathbf{T} [|\varphi u_{\varepsilon}|^2]\right) (\varphi u_{\varepsilon })  \right|^2  + \int_{\mathbb{R}^2} V_{\varepsilon} |\varphi u_{\varepsilon }|^2 \nonumber \\ \nonumber
    & = e_{\varepsilon} + \int_{\mathbb{R}^2} \left| \left( \mathrm{i}\partial_x + \pi \beta f(y) |\varphi(x)|^2 \right) \left( \varphi(x)  u_{\varepsilon}(y)\right) \right|^2 \mathrm{d}x \mathrm{d}y  + \int_{\mathbb{R}} |x|^2 |\varphi(x)|^2 \mathrm{d} x  \\ 
    & = e_{\varepsilon} + \int_{\mathbb{R}^2}    \Big({|\partial_x \varphi |^2 u^2_{\varepsilon}  +  \mathrm{i} \pi \beta (\overline{\varphi} \partial_x \varphi - \varphi\partial_x \overline{\varphi} ) |\varphi|^2 f u^2_{\varepsilon} +  \pi^2 \beta^2|\varphi|^6  f^2u^2_{\varepsilon}  }\Big) + \int_{\mathbb{R}} |x|^2 |\varphi|^2 , \label{integratingY}
\end{align}
where 
\begin{equation}\label{def-f}
    f(y) =  \int_{\mathbb{R}} \operatorname{sgn} (y - \nu ) u^2_{\varepsilon}(\nu) \mathrm{d} \nu . 
\end{equation}
Notice the properties of $f$
$$f' = 2 u_{\varepsilon}^2, \quad f( +\infty) = 1 \quad \text{and} \quad f(-\infty) = -1, $$
so that we obtain
\begin{equation}\label{fProperties}
\begin{split}
    \int_{\mathbb{R}} f u^2 _{\varepsilon} &  =    \frac{1}{2}\int_{\mathbb{R}} f f' =  \frac{1}{4}\int_{\mathbb{R}} (f^2)'  =  \frac{1}{4} (f^2(\infty) - f^2(-\infty))=  0  \\
    \int_{\mathbb{R}} f^2  u^2_{\varepsilon} & =  \frac{1}{2}\int_{\mathbb{R}} f^2 f'  = \frac{1}{6}\int_{\mathbb{R}} (f^3)'   =  \frac{1}{6} (f^3(\infty) - f^3(-\infty))= \frac{1}{3}. 
\end{split}
\end{equation}
Then after integrating out the $y$ variable in \eqref{integratingY}, the 2D energy \eqref{2Denergy} of this trial state $\psi$ becomes
\begin{equation}\nonumber
      \mathcal{E}^{2\mathrm{D}}_{\varepsilon} (\psi) = e_{\varepsilon} +  \mathcal{E}^{1\mathrm{D}}_{} (\varphi),
\end{equation}
where $\mathcal{E}^{1\mathrm{D}}_{}  $ is as in \eqref{1Denergy}, the quintic NLS energy. Running through all possible $\varphi$, we conclude the proof for the energy upper bound stated in Lemma \ref{upperbound}.
\end{proof}

\subsection{Energy lower bound} \label{secLower}

In this section, we prove the following:

\begin{lemma}[Energy lower bound]\label{lowerbound}
Let $E^{2\mathrm{D}}_{\varepsilon} $ and $E^{1\mathrm{D}}_{} $ be as in Definition \ref{GSE}, and let $e_{\varepsilon}$ be as in \eqref{OHenergy}. Then we have
\begin{equation}\nonumber
     \liminf_{\varepsilon \to 0}\LR{{E}^{2\mathrm{D}}_{\varepsilon} - e_{\varepsilon} } \ge  {E}^{1\mathrm{D}}_{}.
\end{equation}
\end{lemma}

One can refer to \cite[Appendix A]{LunRou-15} for the existence of minimizers for the 2D energy \eqref{2Denergy}. 
Let $\Psi$ be a ground state for  \eqref{2Denergy} that we rewrite as
\begin{equation}\label{def-2DPhi}
    \Psi(x,y) = \varphi(x,y) u_{\varepsilon}(y) e^{-\mathrm{i} \beta S[|\varphi|^2 u_\varepsilon^2](x,y)},
\end{equation}
where $u_{\varepsilon}$ is as in \eqref{def-u}, $S[\cdot]$ is as in \eqref{defS} and 
\begin{equation}\label{def-varphi}
    \varphi(x,y) = \Psi(x,y) u_{\varepsilon}^{-1} e^{\mathrm{i} \beta S[|\Psi|^2](x,y)}.
\end{equation}
This function $\varphi$ is well-defined since $u_{\varepsilon}$ vanishes nowhere. 
Then, with the help of a change of gauge \eqref{get-rid-of-phase}, the 2D ground state energy can be written as
\begin{equation}\label{filteredE}
    \begin{split}
        E^{2\mathrm{D}}_{\varepsilon} & = \mathcal{E}^{2\mathrm{D}}_{\varepsilon}(\Psi)  = \int_{\mathbb{R}^2 } \left| \left( - \mathrm{i} \nabla +  \beta  \mathbf{T} [|\varphi |^2 u_{\varepsilon}^2]\right) (\varphi u_{\varepsilon })  \right|^2  + \int_{\mathbb{R}^2} V_{\varepsilon} |\varphi|^2 u_{\varepsilon }^2.
    \end{split}
\end{equation}
Notice that the delta function $\delta_{x=0}$ convoluted 
with $\rho = |\Psi|^2$ is well-defined, because $\sqrt{\rho}$ is in $H^1(\mathbb{R}^2)$ by the diamagnetic inequality \cite[Lemma A.3]{LunRou-15}
\begin{equation}\nonumber
    \int_{\mathbb{R}^2} \abs{\nabla \abs{\Psi}}^2 \le  \int_{\mathbb{R}^2} \abs{\left( - \mathrm{i} \nabla_{} +  \beta \mathbf{A}[|\Psi|^2]\right)\Psi}^2 
\end{equation}
and then $\sqrt{\rho} (x,\cdot)$ is in $H^{\frac{1}{2}}(\mathbb{R})$ by the trace theorem \cite[the 7th Comment on Chapter 9]{Brezis-10}.

\begin{prop}[Energy decoupling]\label{energydecoupling}
    Consider 
    $$ D = - \mathrm{i} \nabla + \mathbf{A} $$
    with real-valued $\mathbf{A}$. Let $\Psi$ be a $L^2$ normalized function such that 
    \begin{equation}\nonumber
        \int_{\mathbb{R}^2} |D \Psi |^2 + \int_{\mathbb{R}^2}   V_{\varepsilon} |\Psi|^2 < \infty
    \end{equation}
    for $V_{\varepsilon}$ as in \eqref{def-v}.
    Then we have
    \begin{equation}\nonumber
        \int_{\mathbb{R}^2} |D \Psi |^2 +V_{\varepsilon} |\Psi|^2  = e_{\varepsilon}  +   \int_{\mathbb{R}^2} |D \varphi |^2 u_{\varepsilon}^2  +   \int_{\mathbb{R}^2}   |x|^2 |\varphi|^2 u_{\varepsilon} ^2,
    \end{equation}
    where $e_{\varepsilon}$ is as in \eqref{OHenergy},  $u_{\varepsilon}$ is as in \eqref{def-u}, and 
    $$ \varphi = \Psi u_{\varepsilon}^{-1}. $$
\end{prop}
\begin{proof}
Similar to the proof of  \cite[Proposition 4.2]{RouYan-23b}.
\end{proof}
Using the energy decoupling in Proposition \ref{energydecoupling}, the energy in \eqref{filteredE} becomes
\begin{equation}\nonumber
    E^{2\mathrm{D}}_{\varepsilon}  = e_{\varepsilon} + \int_{\mathbb{R}^2} \left| \left( - \mathrm{i} \nabla +  \beta \mathbf{T} [|\varphi u_{\varepsilon}|^2]\right) \varphi\right|^2   u_{\varepsilon}^2 + \int_{\mathbb{R}^2} |x|^2 |\varphi|^2 u_{\varepsilon}^2.
\end{equation}
Combining with the energy upper bound in Lemma \ref{upperbound}, we get
\begin{equation}\label{beforeND}
    E^{1\mathrm{D}}_{} \ge E^{2\mathrm{D}}_{\varepsilon} - e_{\varepsilon} = \int_{\mathbb{R}^2} \left| \left( - \mathrm{i} \nabla +  \beta \mathbf{T} [|\varphi u_{\varepsilon}|^2]\right) \varphi\right|^2   u_{\varepsilon}^2 + \int_{\mathbb{R}^2} |x|^2 |\varphi|^2 u_{\varepsilon}^2. 
\end{equation}
We will narrow down the integral domain in the $y$-direction for technical reasons. This will not modify our results due to the harmonic trapping in the $y$-direction. 

Notice that, for a constant ${R} > 0$ ($R$ can depend on $\varepsilon$ and we will let it be like $C \sqrt{\varepsilon}$ later),
\begin{equation}
    \operatorname{curl}{\left(\mathbf{T} [|\varphi u_{\varepsilon}|^2 \mathbbm{1}_{|y|\ge  {R}}]\right)} = \operatorname{curl}{\left(\mathbf{T}_0 * (|\varphi u_{\varepsilon}|^2 \mathbbm{1}_{|y|\ge  {R}})\right)} = -2\pi |\varphi u_{\varepsilon}|^2 \mathbbm{1}_{|y|\ge  {R}} \nonumber
\end{equation}
vanishes on $\{ |y| < {R} \}$. One can check that $\mathbf{T} [|\varphi u_{\varepsilon}|^2 \mathbbm{1}_{|y|\ge  {R}}]$ is in the vector-valued Sobolev space $\mathbf{W}^{1,1}(|y| < {R})$. Hence, according to Theorem \ref{0curl}, there is a function $W_{\varepsilon}$ defined on $\{ |y| < {R} \}$ such that
\begin{equation}\nonumber
    \mathbf{T} [|\varphi u_{\varepsilon}|^2 \mathbbm{1}_{|y|\ge {R}}]
    = \nabla W_{\varepsilon} \quad \text{ on } \quad \{ |y|< {R} \}.
\end{equation}
For convenience, we choose $W_{\varepsilon}$ to be of the form \eqref{defPotential}, i.e.
\begin{equation} \label{def-W}
    W_{\varepsilon}(x,y) \colonequals \int_{\gamma_{(x,y)}} \mathbf{T} [|\varphi u_{\varepsilon}|^2 \mathbbm{1}_{|y|\ge {R}}]
    \cdot \mathrm{d} \mathbf{p},
\end{equation}
where $\gamma_{(x,y)}$ is the piecewise linear curve consisting of the two oriented line segments from the origin $(0,0)$ to $(x,0)$ and from $(x,0)$ to $(x,y)$, and $\mathbf{p}$ denotes the position vector of a point on $\gamma_{(x,y)}$.
Then we have
\begin{equation}\label{property-W}
    W_{\varepsilon}(x,y) = W_{\varepsilon}(x,0) = \pi \int_{0}^x \int_{-\infty}^{\infty} \operatorname{sgn}(\nu) |\varphi(s,\nu)|^2 u^2_{\varepsilon}(\nu) \mathbbm{1}_{|\nu|\ge {R}}(s,\nu) \mathrm{d} \nu \mathrm{d}s.
\end{equation}
Defining $\tilde{\varphi}$ by
\begin{equation}\label{def-tildephi}
    \tilde{\varphi} \colonequals \begin{cases}
        {\varphi} e^{\mathrm{i} \beta W_{\varepsilon}} &\quad \text{ on } \quad \{ |y|<
        {R} \} \\
        \varphi &\quad \text{ on } \quad \{ |y|\ge {R} \} 
    \end{cases} ,
\end{equation}
we then can narrow down the integral in \eqref{beforeND}:
\begin{equation}\nonumber
    \begin{split}
        E^{1\mathrm{D}}_{} & \ge E^{2\mathrm{D}}_{\varepsilon} - e_{\varepsilon} \\
        & \ge  \int_{|y| < R} \left| \left( - \mathrm{i} \nabla +  \beta \mathbf{T}[|\varphi u_{\varepsilon}|^2]\right) \left( \tilde{\varphi} e^{-\mathrm{i} \beta W_{\varepsilon}}\right)\right|^2   u_{\varepsilon}^2 + \int_{\mathbb{R}^2} |x|^2 |\tilde{\varphi}|^2 u_{\varepsilon}^2\\
        & = \int_{|y| < R} \left| \left( - \mathrm{i} \nabla +  \beta \mathbf{T}[|\varphi u_{\varepsilon}|^2] - \beta \nabla W_{\varepsilon}\right) \tilde{\varphi}\right|^2   u_{\varepsilon}^2 + \int_{\mathbb{R}^2} |x|^2 |\tilde{\varphi}|^2 u_{\varepsilon}^2\\
        & = \int_{|y| < R} \left| \left( - \mathrm{i} \nabla +  \beta \mathbf{T}[|\tilde{\varphi} u_{\varepsilon}|^2 \mathbbm{1}_{|y|  < R}]\right) \tilde{\varphi}\right|^2   u_{\varepsilon}^2 + \int_{\mathbb{R}^2} |x|^2 |\tilde{\varphi}|^2 u_{\varepsilon}^2.
    \end{split}
\end{equation}
Making the rescaling $y\to \sqrt{\varepsilon} y$, we have
\begin{equation}\label{afterND}
    E^{1\mathrm{D}}_{} \ge E^{2\mathrm{D}}_{\varepsilon} - e_{\varepsilon} = \int_{|y| < \tilde{R}} \left| \left( \mathrm{i} \partial_x +  \pi \beta f_{\tilde{R},\varepsilon} \right) \tilde{\varphi}_{\varepsilon}\right|^2   u_{1}^2 +  \frac{1}{\varepsilon} \int_{|y| < \tilde{R}} |\partial_y \tilde{\varphi}_{\varepsilon}|^2 u_1^2 +
 \int_{\mathbb{R}^2} |x|^2 |\tilde{\varphi}_{\varepsilon}|^2 u_{1}^2, 
\end{equation}
where
\begin{equation}\label{def-tildephiepsilon}
\tilde{\varphi}_{\varepsilon} (x,y) = \tilde{\varphi} \left(x, \sqrt{\varepsilon} y\right)
\end{equation}
and
\begin{equation}\label{def-fR}
    \tilde{R} = \frac{R}{\sqrt{\varepsilon}}, \quad
    f_{\tilde{R},\varepsilon} (x,y) = \int_{-\tilde{R}}^{\tilde{R}}\operatorname{sgn}(y -\nu) |\tilde{\varphi}_{\varepsilon}(x,\nu)|^2 u_1^2(\nu) \mathrm{d} \nu.
\end{equation}
We would like to fix $\tilde{R}$ as a constant independent of $\varepsilon$, so we set $$R = \tilde{R} \sqrt{\varepsilon}.$$
With direct calculations, we find the following properties for $f_{\tilde{R},\varepsilon}$,
\begin{equation}\label{fRproperty}
    \partial_y f_{\tilde{R},\varepsilon} = 2|\tilde{\varphi}_{\varepsilon}|^2
    u_{1}^2 \mathbbm{1}_{|y| < \tilde{R}}, \quad f_{\tilde{R},\varepsilon}(x,\infty) =\int_{-\tilde{R}}^{\tilde{R}}|\tilde{\varphi}_{\varepsilon}(x,\nu)|^2 u_1^2(\nu) \mathrm{d} \nu = - f_{\tilde{R},\varepsilon}(x,-\infty).
\end{equation}

Applying the diamagnetic inequality \cite[Theorem 2.1.1]{FourHel-10}
\begin{equation}\nonumber
    \left| \partial_x |\tilde{\varphi}_{\varepsilon}|\right|^2 \le \left| \left( \mathrm{i} \partial_x +  \pi \beta f_{\tilde{R},\varepsilon} \right) \tilde{\varphi}_{\varepsilon}\right|^2
\end{equation}
to the first square term on the right-hand side of \eqref{afterND}, we obtain
\begin{equation}\label{controlonnabla}
    \int_{|y| < \tilde{R}}\left| \nabla |\tilde{\varphi}_{\varepsilon}|\right|^2 u_1^2 \le E^{1\mathrm{D}}_{}.
\end{equation}
In order to pass the limit, we still need a similar control of $\nabla \tilde{\varphi}_{\varepsilon}$.
Expanding the first square term on the right-hand side of \eqref{afterND}
\begin{multline}\label{writeout}
     \int_{|y| < \tilde{R}} \left| \left( \mathrm{i} \partial_x +  \pi \beta f_{\tilde{R},\varepsilon} \right) \tilde{\varphi}_{\varepsilon}\right|^2   u_{1}^2 = \int_{|y| < \tilde{R}} |\partial_x \tilde{\varphi}_{\varepsilon}|^2 u_1^2 \\ + \mathrm{i} \pi \beta \int_{|y| < \tilde{R}} f_{\tilde{R} , \varepsilon} \left( \overline{\tilde{\varphi}_{\varepsilon}} \partial_x \tilde{\varphi}_{\varepsilon} - \tilde{\varphi}_{\varepsilon} \partial_x \overline{\tilde{\varphi}_{\varepsilon}}  \right)u_1^2   + \pi^2 \beta^2 \int_{|y| < \tilde{R}} f_{\tilde{R} , \varepsilon}^2|\tilde{\varphi}_{\varepsilon}|^2 u_1^2, 
\end{multline}
and applying the inequality 
$$ 2 |ab| \le |a|^2 +|b|^2  $$ to the cross-term, we find, for a constant $\eta \in (0,1)$ independent of $\varepsilon$,
\begin{multline}
    \int_{|y| < \tilde{R}} \left| \left( \mathrm{i} \partial_x +  \pi \beta f_{\tilde{R},\varepsilon} \right) \tilde{\varphi}_{\varepsilon}\right|^2   u_{1}^2 \\ \ge  \left( 1-\eta^2 \right) \int_{|y| < \tilde{R}} |\partial_x \tilde{\varphi}_{\varepsilon}|^2 u_1^2 - \left( \frac{1}{\eta^2} -1 \right) \pi^2 \beta^2 \int_{|y| < \tilde{R}} f_{\tilde{R} , \varepsilon}^2|\tilde{\varphi}_{\varepsilon}|^2 u_1^2. \label{afterAMGM}
\end{multline}
To obtain a nice control on the first term, i.e. the term related to $\partial_x \tilde{\varphi}_{\varepsilon}$, we claim the following control on the last square term on the right-hand side of \eqref{afterAMGM}:

\begin{lemma}\label{fRcontrol}
    With the notation above, we have
    \begin{equation}\nonumber
        \int_{|y| < \tilde{R}} f_{\tilde{R} , \varepsilon}^2|\tilde{\varphi}_{\varepsilon}|^2 u_1^2 \le C
    \end{equation}
    for a constant $C$ independent of $\varepsilon$.
\end{lemma}
\begin{proof}
Recall the properties of $f_{\tilde{R} , \varepsilon}$ in \eqref{fRproperty},
which lead to
\begin{multline}\nonumber
        \int_{|y| < \tilde{R}} f_{\tilde{R} , \varepsilon}^2|\tilde{\varphi}_{\varepsilon}|^2 u_1^2  = \frac{1}{2} \int_{\mathbb{R}^2} f_{\tilde{R} , \varepsilon}^2 \partial_y f_{\tilde{R} , \varepsilon} 
        = \frac{1}{6}\int_{\mathbb{R}^2} \partial_y f_{\tilde{R} , \varepsilon}^3
        = \frac{1}{6} \int_{\mathbb{R}}\left[ f_{\tilde{R} , \varepsilon}^3(x,\infty) - f_{\tilde{R} , \varepsilon}^3(x,-\infty)\right]\mathrm{d} x \\
        =  \frac{C^3_{\tilde{R}}}{3} \int_{\mathbb{R}} \left(\int_{-\tilde{R}}^{\tilde{R}}|\tilde{\varphi}_{\varepsilon}(x,\nu)|^2 \frac{u_1^2(\nu)}{C_{\tilde{R}}} \mathrm{d} \nu \right)^3\mathrm{d} x 
\end{multline}
with 
\begin{equation} \nonumber
    C_{\tilde{R}} = \int_{-\tilde{R}}^{\tilde{R}} u_1^2(\nu) \mathrm{d} \nu  \le 1.
\end{equation}
Using Jensen's inequality, we obtain 
\begin{equation}\nonumber
    \begin{split}
        \int_{|y| < \tilde{R}} f_{\tilde{R} , \varepsilon}^2|\tilde{\varphi}_{\varepsilon}|^2 u_1^2 & \le \frac{C_{\tilde{R}}^2}{3} \int_{\mathbb{R}} \int_{-\tilde{R}}^{\tilde{R}}|\tilde{\varphi}_{\varepsilon}(x,\nu)|^6 u_1^2(\nu) \mathrm{d} \nu \mathrm{d} x.
    \end{split}
\end{equation}
Then thanks to the Sobolev embedding $H^{1}(\Omega\subset\mathbb{R}^2) \subset L^6(\Omega) $ \cite[Corollary 9.14]{Brezis-10} and the definition of $u_{1}$ as in \eqref{def-u}, we can control the square term further:
\begin{multline}\nonumber
        \int_{|y| < \tilde{R}} f_{\tilde{R} , \varepsilon}^2|\tilde{\varphi}_{\varepsilon}|^2 u_1^2  \le \int_{|y| < \tilde{R}} |\tilde{\varphi}_{\varepsilon}|^6 
        \le C \cdot \left(\int_{|y| < \tilde{R}}\left| \nabla |\tilde{\varphi}_{\varepsilon}|\right|^2 + \int_{|y| < \tilde{R}}   |\tilde{\varphi}_{\varepsilon}|^2 \right)^3 \\
        \le  \frac{C}{u_1^6\left(\tilde{R}\right)}  \left(\int_{|y| < \tilde{R}}\left| \nabla |\tilde{\varphi}_{\varepsilon}|\right|^2 u_1^2 + \int_{|y| < \tilde{R}}   |\tilde{\varphi}_{\varepsilon}|^2 u_1^2 \right)^3,
\end{multline}
where $C$ is a constant from the Sobolev inequality and only depends on $\tilde{R}$. 
Recalling the control on $\nabla |\tilde{\varphi}_{\varepsilon}|$ in \eqref{controlonnabla}, we find
\begin{equation}
    \int_{|y| < \tilde{R}} f_{\tilde{R} , \varepsilon}^2|\tilde{\varphi}_{\varepsilon}|^2 u_1^2 \le  \frac{C}{u_1^6({\tilde{R}})}  \left( E^{1\mathrm{D}}_{} +1\right)^3,
\label{control-fsq}
\end{equation}
which completes the proof of Lemma \ref{fRcontrol}.
\end{proof}

We denote by
$$\Omega_{\tilde{R}} = \mathbb{R}\times  (- \tilde{R}, \tilde{R} ),$$ 
and we define the weighted Sobolev space as following:
\begin{equation}\nonumber
    H_{u_1}^1(\Omega_{\tilde{R}}) = H_{}^1\left(\Omega_{\tilde{R}};u_1^2(y) \mathrm{d}x \mathrm{d} y\right), \quad L^p_{u_1} (\Omega_{\tilde{R}}) = L^p\left(\Omega_{\tilde{R}} ;u_1^2(y) \mathrm{d}x \mathrm{d} y\right).
\end{equation}

According to the same philosophy as \cite[Lemma 12]{Lewin-compact}, we have the following lemma:
\begin{lemma}[Extracting the locally convergent part] \label{locallycvg}
    Let $(\varphi_n )_n$ be a sequence in $H^1_{u_1}(\Omega_{\tilde{R}})$ such that $\varphi_n \rightharpoonup \varphi$ weakly in $H^1_{u_1}(\Omega_{\tilde{R}})$ and let $M_k \ge 0$ such that $M_k \to \infty$ as $k\to \infty$. Then there exists a subsequence $( \varphi_{n_{k}})_{k}$ such that
    \begin{equation}\nonumber
        \int_{\Omega_{\tilde{R}} \cap\{|{{x}}| \le M_k \}} |\varphi_{n_k}|^2 u_1^2\to \int_{\Omega_{\tilde{R}}} |\varphi|^2 u_1^2 
    \end{equation}
    as $k\to \infty$. In particular, we have $\varphi_{n_k} \mathbbm{1}_{|{{x}}| \le M_k } \to \varphi$ strongly in $L^p_{u_1} (\Omega_{\tilde{R}})$ for all $2 \le p < \infty$.
\end{lemma}

\begin{prop}[Convergence]\label{cvgprop}
    With the notation above, there exists $\varphi_0  \in H^1(\mathbb{R})$ and a vanishing sequence $(\varepsilon_k)_{k\in \mathbb{N}}$ both independent of $\tilde{R}$ such that
\begin{equation}\nonumber
    \tilde{\varphi}_{\varepsilon_k} \xrightarrow[k\to \infty]{} {\varphi}_{0} \quad \text{weakly in} \quad H_{u_1}^1 (\Omega_{\tilde{R} }),
\end{equation} 
and
\begin{equation}\label{rescaled-lpcvg}
    \tilde{\varphi}_{\varepsilon_k} \xrightarrow[k\to \infty]{} \varphi_{0} \quad \text{strongly in} \quad L^p_{u_1} (\Omega_{\tilde{R} })  \quad \text{for} \quad p\ge 2.
\end{equation}
\end{prop}

\begin{proof}
Inequalities \eqref{afterND}, \eqref{afterAMGM} and \eqref{control-fsq} result in the following estimate
\begin{equation}\nonumber
    \int_{|y| < { \tilde{R} }} |\partial_x \tilde{\varphi}_{\varepsilon}|^2 u_1^2 \le \frac{1}{1-\eta^2} \left( E^{1\mathrm{D}}_{\beta } + \frac{(1-\eta^2) C}{\eta^2 u_1^6({{ \tilde{R} }})}  \pi^2 \beta^2 \left( E^{1\mathrm{D}}_{} +1\right)^3 \right)
\end{equation}
for any constant $ \tilde{R}  >0$ and $\eta \in (0,1)$. Again from the estimate as in \eqref{afterND}, we know 
\begin{equation}
    \int_{|y|<  \tilde{R} } |\partial_y \tilde{\varphi}_{\varepsilon}|^2 u_1^2 \le E^{1\mathrm{D}}_{} \varepsilon. \label{partialY}
\end{equation}
Hence, $\{\tilde{\varphi}_{\varepsilon}\}_{\varepsilon}$ is bounded in the weighted Sobolev space $H_{u_1}^1(\Omega_{\tilde{R} })$.
Then, there exists a vanishing sequence $(\varepsilon_k)_{k\in \mathbb{N}}$ such that $\tilde{\varphi}_{\varepsilon_k} $ converges weakly in $H_{u_1}^1 (\Omega_{\tilde{R} })$ to some $ {\varphi}_{\tilde{R} } $ as $k$ goes to $\infty$.
The limit $ {\varphi}_{\tilde{R} } $ only depends on $x$ due to the control in \eqref{partialY}.
So this shows that
\begin{equation}\label{deri}
    \liminf_{k \to \infty} \int_{|y| < \tilde{R} } |\partial_x \tilde{\varphi}_{\varepsilon_k}|^2 u_1^2  \ge  \int_{\mathbb{R}} |\partial_x {\varphi}_{\tilde{R} }|^2 \cdot \int_{-\tilde{R} }^{\tilde{R} } u_1^2,
\end{equation}
which implies that ${\varphi}_{\tilde{R} } \in H^1(\mathbb{R})$.

With the help of Lemma \ref{locallycvg} and the same sequence $(M_k)_{k\in\mathbb{N}}$, after passing to a subsequence, we have
\begin{equation}\nonumber
    \tilde{\varphi}_{\varepsilon_k} \mathbbm{1}_{|x|\le M_k} \xrightarrow[ k\to \infty]{} \varphi_{\tilde{R} } \quad \text{strongly in} \quad L^p_{u_1} (\Omega_{\tilde{R} }) \quad\text{for} \quad p\ge 2.
\end{equation}   
On the other hand, 
\begin{equation}\nonumber
    \begin{split}
        \int_{\Omega_{\tilde{R} }} | \tilde{\varphi}_{\varepsilon_k} \mathbbm{1}_{|x|> M_k}|^p u_1^2  
        & \le  \frac{1}{M_k}  \int_{\Omega_{\tilde{R} } \cap\{|{{x}}| > M_k \}}  |x| | \tilde{\varphi}_{\varepsilon_k}|^p u_1^2  \\
        &  \le  \frac{1}{M_k}  \left( \int_{\Omega_{\tilde{R} } \cap\{|{{x}}| \le M_k \}}   |x|^2 | \tilde{\varphi}_{\varepsilon_k}|^2 u_1^2 \right)^{\frac{1}{2}} \left( \int_{\Omega_{\tilde{R} } \cap\{|{{x}}| \le M_k \}}  | \tilde{\varphi}_{\varepsilon_k}|^{2(p-1)}u_1^2 \right)^{\frac{1}{2}} 
    \end{split}
\end{equation}
vanishes as $k\to \infty$ due to the control in \eqref{afterND} and Sobolev embeddings $$ H^1_{u_1}(\Omega_{\tilde{R} }) \subset L^{2(p-1)}_{u_1} (\Omega_{\tilde{R} }) \quad \text{for} \quad p\ge 2.$$ 
Therefore, $\tilde{\varphi}_{\varepsilon_k}$ converges strongly in $L^p_{u_1} (\Omega_{\tilde{R} })  $ to $\varphi_{\tilde{R} } $ for $p\ge 2$ as $k$ goes to $\infty$.

Recall the definition of $W_{\varepsilon}$ in \eqref{def-W} that depends on $R$, i.e. different value of $R$ corresponds to different $W_{\varepsilon}$.
To clarify it for the moment, we denote $W_{\varepsilon}$, $\tilde{\varphi}_{}$ and $\tilde{\varphi}_{\varepsilon} $ by $W^{(R)}_{\varepsilon}$, $\tilde{\varphi}^{(R)}$ and $\tilde{\varphi}_{\varepsilon}^{(\tilde{R})} $ respectively, where ${R} =  \tilde{R}\sqrt{\varepsilon}$ as above.

Let $\tilde{R}_n >0$ such that $\tilde{R}_n \to \infty$ as $n\to \infty$ and $\tilde{R}_n < \tilde{R}_{n+1}$ for $n\in \mathbb{N}$. 
Using the diagonal method, after passing to a subsequence, $(\tilde{\varphi}_{\varepsilon_k}^{(\tilde{R}_n)})_{k\in\mathbb{N}}$ converges strongly in $L^p_{u_1} (\Omega_{\tilde{R}_n})  $ to $\varphi_{\tilde{R}_n}$ for  all $n \in \mathbb{N}$ with a same vanishing sequence $(\varepsilon_k)_{k\in \mathbb{N}}$, and it also converges pointwise almost everywhere to the limit $\varphi_{\tilde{R}_n}$ that is independent of $y$.
For any constant $\tilde{R} >0$, there exists $n\in \mathbb{N}$ such that $\tilde{R} < \tilde{R}_n$.
Recall the expression of $W_{\varepsilon}^{(R)}$ in \eqref{property-W}.
We have
\begin{equation}\nonumber
    \begin{split}
        W_{\varepsilon}^{({R}_n)} - W_{\varepsilon}^{({R} )}
        & = - \pi {\int_{0}^x \int_{-\infty}^{\infty} \operatorname{sgn}(\nu) { |{\varphi}_{}(s,\nu)|^2} u^2_{\varepsilon}(\nu) \mathbbm{1}_{R \le|\nu|\le {R}_n }(s,\nu) \mathrm{d} \nu \mathrm{d}s}\\
        & = - \pi {\int_{0}^x \int_{-\infty}^{\infty} \operatorname{sgn}(\nu) { |\tilde{\varphi}_{\varepsilon}^{(\tilde{R}_n )}(s,\nu)|^2} u^2_{1}(\nu) \mathbbm{1}_{\tilde{R} \le|\nu|\le \tilde{R}_n }(s,\nu) \mathrm{d} \nu \mathrm{d}s}.
    \end{split}
\end{equation}
Since $(\tilde{\varphi}_{\varepsilon_k}^{(\tilde{R}_n)})_{k\in\mathbb{N}}$ converges strongly in $L^2_{u_1} (\Omega_{\tilde{R}_n })$ to $\varphi_{\tilde{R}_n }$, we obtain
\begin{equation}\label{Wdiff}
    \lim_{k\to \infty} \LR{ W_{\varepsilon_k}^{({R}_n)} - W_{\varepsilon_k}^{({R} )} } = - \pi {\int_{0}^x \int_{-\infty}^{\infty} \operatorname{sgn}(\nu) { |{\varphi}_{\tilde{R} }(s)|^2} u^2_{1}(\nu) \mathbbm{1}_{\tilde{R} \le|\nu|\le \tilde{R}_n }(s,\nu) \mathrm{d} \nu \mathrm{d}s} = 0,
\end{equation}
i.e. $  W_{\varepsilon_k}^{({R}_n)} - W_{\varepsilon_k}^{({R} )} $ vanishes pointwise as $k\to 0$.
With direct calculations, it turns out that
\begin{equation}\nonumber
\begin{split}
    \norm{ \tilde{\varphi}_{\varepsilon_k}^{(\tilde{R})} - \varphi_{\tilde{R}_n}}_{L^p_{u_1}(\Omega_{\tilde{R}})}  & = \norm{\tilde{\varphi}_{\varepsilon_k}^{(\tilde{R}_n)} e^{-\mathrm{i} \beta \LR{W_{\varepsilon_k}^{(R_n)} - W_{\varepsilon_k}^{(R)}}} -  \varphi_{\tilde{R}_n} }_{L^p_{u_1}(\Omega_{\tilde{R}})} \\
    & \le  \norm{\tilde{\varphi}_{\varepsilon_k}^{(\tilde{R}_n)} - \varphi_{\tilde{R}_n}}_{L^p_{u_1}(\Omega_{\tilde{R}})} +  \norm{\varphi_{\tilde{R}_n}\LR{e^{- \mathrm{i}\beta \LR{W_{\varepsilon_k}^{({R}_n)} - W_{\varepsilon_k}^{({R})}}}-1 }}_{L^p_{u_1}(\Omega_{\tilde{R}})}
\end{split}
\end{equation}
vanishes as $k\to \infty$, where we apply the dominated convergence theorem to the last term on the right-hand side.
Thus, $(\tilde{\varphi}_{\varepsilon_k}^{(\tilde{R})} )_{k\in \mathbb{N}}$ converges to $\varphi_{\tilde{R}_n}$ in $L^2_{u_1}(\Omega_{\tilde{R}})$ for all $\tilde{R} <\tilde{R}_n$, which implies that 
$$ \varphi_{\tilde{R}_n} = \varphi_{\tilde{R}_m} $$
for all $n, m \in \mathbb{N}$, i.e. $\varphi_{\tilde{R}_n}$ does not depend on $\tilde{R}_n$ and we denote it by $\varphi_0$ in the sequel.
\end{proof}

\begin{cor}
    With the same notation above, consider $f_{\tilde{R}, \varepsilon}$ defined in \eqref{def-fR}. Then we have
\begin{equation}\label{fcvg}
    f_{\tilde{R},\varepsilon_k} \xrightarrow[k\to \infty]{} f_{\tilde{R}} \quad \text{strongly in} \quad L^q_{u_1}(\Omega_{\tilde{R}}) \quad \text{for} \quad q\ge 1,
\end{equation}
where 
\begin{equation}\nonumber
    f_{\tilde{R}} (x,y) \colonequals  |{\varphi}_{0}(x)|^2  \int_{-\tilde{R}}^{\tilde{R}}\operatorname{sgn}(y -\nu)  u_1^2(\nu) \mathrm{d} \nu.
\end{equation}
\end{cor}
\begin{proof}
Using similar properties as in \eqref{fRproperty}, we obtain
\begin{multline}\nonumber
    \int_{-\tilde{R}}^{\tilde{R}} f_{\tilde{R}}^2(x,y)  u_1^2 (y) \mathrm{d} y = \frac{1}{2|\varphi_0(x)|^2} \int_{\mathbb{R}} f_{\tilde{R}}^2(x,y) \partial_y f_{\tilde{R}} (x,y)  \mathrm{d} y = \frac{1}{6|\varphi_0(x)|^2} \int_{\mathbb{R}} \partial_y f_{\tilde{R}}^3 (x,y)  \mathrm{d} y \\ = \frac{f_{\tilde{R}}^3 (x,\infty) - f_{\tilde{R}}^3 (x,-\infty) }{6|\varphi_0(x)|^2}  = \frac{1}{3} |\varphi_0(x)|^4 \left(\int_{-\tilde{R}}^{\tilde{R}} u_1^2  \right)^3
\end{multline}
and
\begin{multline}\nonumber
     \int_{-\tilde{R}}^{\tilde{R}} f_{\tilde{R}}(x,y)  u_1^2 (y) \mathrm{d} y = \frac{1}{2|\varphi_0(x)|^2} \int_{\mathbb{R}} f_{\tilde{R}}(x,y) \partial_y f_{\tilde{R}} (x,y) \mathrm{d} y = \frac{1}{4|\varphi_0(x)|^2} \int_{\mathbb{R}} \partial_y f_{\tilde{R}}^2 (x,y)  \mathrm{d} y \\ = \frac{f_{\tilde{R}}^2 (x,\infty) - f_{\tilde{R}}^2 (x,-\infty) }{4|\varphi_0(x)|^2}  = 0.
\end{multline}
The convergence in \eqref{fcvg} results from the following calculations:
\begin{equation}\nonumber
    \begin{split}
       \int_{|y|\le \tilde{R}} |f_{\tilde{R},\varepsilon_k} - f_{\tilde{R}} |^q u_1^2 = &  \int_{|y|\le \tilde{R}} \left| \int_{-\tilde{R}}^{\tilde{R}}\operatorname{sgn}(y -\nu)\left( |\tilde{\varphi}_{\varepsilon_k}(x,\nu)|^2 - |{\varphi}_{0}(x)|^2 \right) u_1^2(\nu) \mathrm{d} \nu \right|^q u_1^2(y) \mathrm{d}x \mathrm{d}y\\
        \le & \int_{|y|\le \tilde{R}} \left( \int_{-\tilde{R}}^{\tilde{R}}\left| |\tilde{\varphi}_{\varepsilon_k}(x,\nu)|^2 - |{\varphi}_{0}(x)|^2 \right| u_1^2(\nu) \mathrm{d} \nu \right)^q u_1^2(y) \mathrm{d}x \mathrm{d}y\\
        \text{(Jensen's inequality)} \quad\le &\int_{|y|\le \tilde{R}} \left| |\tilde{\varphi}_{\varepsilon_k}(x,\nu)|^2 - |{\varphi}_{0}(x)|^2 \right|^q u_1^2(\nu) \mathrm{d} \nu \mathrm{d}x\\
        \text{(H\"older's inequality)} \quad\le &  \left( \int_{|y|\le \tilde{R}}  \left( |\tilde{\varphi}_{\varepsilon_k}| + |{\varphi}_{0}| \right)^{2q} u_1^2  \right)^{\frac{1}{2}}  \left( \int_{|y|\le \tilde{R}} \left| \tilde{\varphi}_{\varepsilon_k}
        - {\varphi}_{0} \right|^{2q} u_1^2  \right)^{\frac{1}{2} } \\
        \le  &  2^{q -1} \left( \norm{\tilde{\varphi}_{\varepsilon_k}}_{L^{2q}_{u_1}(\Omega_{\tilde{R}})}^q  + \norm{{\varphi}_{0}}^q_{L^{2q}_{u_1}(\Omega_{\tilde{R}})}  \right)  \norm{ \tilde{\varphi}_{\varepsilon_k}
        - {\varphi}_{0} }^q_{L^{2q}_{u_1}(\Omega_{\tilde{R}})} \xrightarrow[k\to \infty]{} 0
    \end{split}
\end{equation}
due to the strong convergence of $\tilde{\varphi}_{\varepsilon_k}$ in $L^{2q}_{u_1}(\Omega_{\tilde{R}})$ stated in Proposition \ref{cvgprop}.    
\end{proof}

Thanks to the strong convergences of $\tilde{\varphi}_{\varepsilon_k}$ and $f_{\tilde{R},\varepsilon_k}$ stated in \eqref{rescaled-lpcvg} and \eqref{fcvg}, we can quickly obtain that
\begin{equation}\nonumber
    f_{\tilde{R},\varepsilon_k}\tilde{\varphi}_{\varepsilon_k} \xrightarrow[k \to \infty]{} f_{\tilde{R}}\varphi_0 \quad \text{strongly in} \quad L^2_{u_1}(\Omega_{\tilde{R}}),
\end{equation} which implies that 
\begin{equation}\label{f^2}
    \lim_k \int_{|y|\le\tilde{R}} f_{\tilde{R} , \varepsilon_k}^2|\tilde{\varphi}_{\varepsilon_k}|^2 u_1^2 =  \int_{|y|\le \tilde{R}} f_{\tilde{R}}^2|{\varphi}_{0}|^2 u_1^2 =  \frac{1}{3} \left( \int_{\mathbb{R}} |{\varphi}_{0}|^6 \right)\cdot \left(\int_{-\tilde{R}}^{\tilde{R}} u_1^2 \right)^3,
\end{equation}
\begin{equation}\label{crossterm}
    \lim_k \int_{|y|\le \tilde{R}} f_{\tilde{R} , \varepsilon_k} \left( \overline{\tilde{\varphi}_{\varepsilon_k}} \partial_x \tilde{\varphi}_{\varepsilon_k} - \tilde{\varphi}_{\varepsilon_k} \partial_x \overline{\tilde{\varphi}_{\varepsilon_k}}  \right)u_1^2 = \int_{|y|\le \tilde{R}} f_{\tilde{R}} \left( \overline{{\varphi}_{0}} \partial_x {\varphi}_{0} - {\varphi}_{0} \partial_x \overline{{\varphi}_{0}}  \right)u_1^2 = 0.
\end{equation}
Applying Fatou's lemma, we have
\begin{equation}\label{potential}
   \liminf_{k \to \infty} \int_{\mathbb{R}^2} |x|^2 |\tilde{\varphi}_{\varepsilon_k}|^2 u_{1}^2 \ge \int_{\mathbb{R}} |x|^2 |{\varphi}_{0}|^2.
\end{equation}
In conclusion, it results from \eqref{afterND},\eqref{writeout}
,\eqref{deri},  \eqref{f^2}, \eqref{crossterm} and \eqref{potential} that
\begin{equation}\nonumber
   \liminf_{\varepsilon \to 0} \left(E^{2\mathrm{D}}_{\varepsilon} -e_{\varepsilon} \right)\ge \int_{\mathbb{R}} |\partial_x {\varphi}_{0}|^2 \cdot \int_{-\tilde{R}}^{\tilde{R}} u_1^2 +\frac{1}{3} \pi^2 \beta^2 \int_{\mathbb{R}} |{\varphi}_{0}|^6 \cdot \left(\int_{-\tilde{R}}^{\tilde{R}} u_1^2\right)^3 + \int_{\mathbb{R}} |x|^2 |{\varphi}_{0}|^2 
\end{equation}
for any $\tilde{R}>0$. 
Letting $\tilde{R}$ go to infinity, we find
\begin{equation}\label{conclu}
     \liminf_{\varepsilon \to 0} \left(E^{2\mathrm{D}}_{\varepsilon} -e_{\varepsilon} \right)\ge \int_{\mathbb{R}} |\partial_x {\varphi}_{0}|^2+\frac{1}{3} \pi^2 \beta^2 \int_{\mathbb{R}} |{\varphi}_{0}|^6  + \int_{\mathbb{R}} |x|^2 |{\varphi}_{0}|^2 \ge E^{1\mathrm{D}}_{}.
\end{equation}
The last inequality of \eqref{conclu} requires the limit ${\varphi}_{0}$ to be normalized in $L^2$, and it is indeed normalized due to a similar proof as in \cite[Lemma 4.10]{RouYan-23b}. This concludes the proof for the energy lower bound stated in Lemma \ref{lowerbound}.

\subsection{Convergence of ground states}   \label{SecL2cvg}

Using the same notation as in the previous section, we now complete the proof of Theorem \ref{gsfunction}.

\begin{lemma}[Vanishing phase factor] \label{vanishW}
    Let $W_{\varepsilon}$ be defined as in \eqref{def-W} and consider the sequence $(\varepsilon_k)_{k\in \mathbb{N}}$ in Proposition \ref{cvgprop}. Then $ W_{\varepsilon_k}$ vanishes pointwise as $k$ goes to 0.
\end{lemma}
\begin{proof}
To avoid the ambiguity, we denote $W_{\varepsilon}$ by $W^{(R)}_{\varepsilon}$ in this proof.

Recall the expression of $W_{\varepsilon}^{(R)}$ in \eqref{property-W}. 
Applying the energy upper bound in Lemma \ref{upperbound} and the energy equation in \eqref{filteredE}, we have
\begin{equation}\nonumber
    \abs{W_{\varepsilon}^{(R)} }\le \pi \int_{|y| \ge \tilde{R} \sqrt{\varepsilon}} |\varphi|^2 u^2_{\varepsilon} \le \frac{\pi}{\tilde{R}^2 {\varepsilon}} \int_{|y| \ge  \tilde{R} \sqrt{\varepsilon}} |y|^2 |\varphi|^2 u^2_{\varepsilon} \le \frac{C}{\tilde{R}^2}
\end{equation}
for some constant $C$ independent of $\varepsilon$ and $\tilde{R}$.
Using \eqref{Wdiff} in the proof of Proposition \ref{cvgprop}, it turns out that
\begin{equation}\nonumber
    \limsup_{k\to \infty}\abs{W_{\varepsilon_k}^{(R)}} \le \limsup_{k\to \infty}\abs{W_{\varepsilon_k}^{(R_n)}} + \lim_{k\to \infty}\abs{W_{\varepsilon_k}^{(R_n)} - W_{\varepsilon_k}^{(R)}} \le \frac{C}{\tilde{R}_n}
\end{equation}
for any $n\in \mathbb{N}$.
Therefore, $ W_{\varepsilon_k}$ vanishes pointwise as $k$ goes to 0.
\end{proof}

The vanishing phase factor $W_{\varepsilon_k}$ implies that the phase $e^{\mathrm{i}\beta W_{\varepsilon_k}}$ has no effect on the $L^p_{u_1}$ limit of $(\tilde{\varphi}^{}_{\varepsilon_k})_{k\in \mathbb{N}}$. 
Thus, $({\varphi}_{\varepsilon_k})_{k\in \mathbb{N}}$ also converges strongly to $\varphi_0$ in $L^p_{u_1}\LR{\Omega_{\tilde{R}}}  $ for any constant $\tilde{R}>0$, where $\varphi$ defined as in \eqref{def-varphi} is independent of $\tilde{R}$ and
\begin{equation}\nonumber
    \varphi_{\varepsilon}(x,y) = \varphi \LR{x, \sqrt{\varepsilon} y}.
\end{equation}
We claim that $({\varphi}_{\varepsilon_k})_{k\in \mathbb{N}}$ converges strongly to $\varphi_0$ in $L^2_{u_1}(\mathbb{R}^2)$. 
It immediately follows from the fact that 
\begin{equation}\label{cvgWhole}
    \int_{|y|\ge \tilde{R}} \abs{{\varphi}_{\varepsilon}}^2 u_1^2 = \int_{|y|\ge \tilde{R} \sqrt{\varepsilon} } \abs{{\varphi}_{}}^2 u_{\varepsilon}^2 \le \frac{1}{\tilde{R}^2 \varepsilon} \int_{|y|\ge \tilde{R} \sqrt{\varepsilon} } |y|^2 \abs{{\varphi}_{}}^2 u_{\varepsilon}^2  \le \frac{C}{\tilde{R}^2}
\end{equation} 
for some constant $C$ independent of $\varepsilon$ and $\tilde{R}$, where the energy upper bound in Lemma \ref{upperbound} and the energy equation in \eqref{filteredE} are used to obtain the right-hand side of \eqref{cvgWhole}.
Recalling the definitions in \eqref{def-Sepsilon} and \eqref{def-2DPhi}, we have that
\begin{equation}\nonumber
    \begin{split}
        \norm{\Psi - \varphi_0 u_{\varepsilon_k}}_{L^2(\mathbb{R}^2)} & = \norm{{\varphi} u_{\varepsilon_k}  e^{-\mathrm{i}\beta {S\left[|{\varphi}|^2u^2_{\varepsilon_k}\right]}} - \varphi_0 u_{\varepsilon_k}}_{L^2(\mathbb{R}^2)}\\
        & \le \norm{ {\varphi}_{\varepsilon_k} - \varphi_0}_{L^2_{u_1}(\mathbb{R}^2)} + \norm{\varphi_0  \LR{e^{-\mathrm{i}\beta {S_{\varepsilon_k}\left[| {\varphi}_{\varepsilon_k}|^2u^2_{1}\right]}}  -1 }}_{L^2_{u_1}(\mathbb{R}^2)}
    \end{split}
\end{equation}
vanishes as $k$ goes to 0, where the phase factor $S_{\varepsilon_k}[|{\varphi}_{\varepsilon_k}|^2u^2_{1}]$ also vanishes pointwise and we again apply the dominated convergence theorem to the last term on the right-hand side.
Therefore, we have completed the proof for Theorem \ref{gsfunction}.

\section{Dimensional reduction for time-dependent solutions} 

In this section, as stated in Theorem \ref{dynamicresult}, we will prove that the solution to the initial value problem \eqref{IVP-A} decomposes into the form of $e^{-i\tilde{S}}\varphi u_{\varepsilon}$ as $\varepsilon \to 0$, where $\varphi$ is the solution to the initial value problem \eqref{1Divp}. 
In other words, the dimensional reduction by the strong confinement on the Chern-Simons-Schr\"odinger equation yields the quintic nonlinear Schr\"odinger equation as the 1D limit model. 
We first provide heuristic calculations to give a flavor of why the limiting equation is the quintic NLS. 
Then we present a detailed proof.

\subsection{Heuristic calculations} \label{trialPDE}
Assume the solution to the 2D anyonic equation \eqref{2DPDE} is in the form 
\begin{equation}\nonumber
    \psi(t,x,y) = \varphi(t,x) u_{\varepsilon}(y) e^{-\mathrm{i} e_{\varepsilon} t - \mathrm{i} \beta S[|\varphi u_{\varepsilon} |^2](t,x,y)}.
\end{equation}
In this section, we compute the equation that $\varphi$ must satisfy, namely, the quintic defocusing NLS equation 
\begin{equation}\nonumber
        \mathrm{i} \partial_t \varphi = -\partial_x^2 \varphi + |x|^2 \varphi + \pi^2 \beta^2 |\varphi|^4 \varphi .
\end{equation}
The main idea is to multiply both sides of the 2D equation \eqref{2DPDE} by 
$${u_{\varepsilon}}e^{\mathrm{i} e_{\varepsilon} t + \mathrm{i} \beta S[|\varphi u_{\varepsilon} |^2]}$$ 
and then integrate over $y$.
The calculations are as follows.\\

\noindent$\bullet$ {\bf{The time derivative part.}}
We first trace out the effect in the tight confinement direction:
\begin{equation}\nonumber
\begin{split}
    \int_{\mathbb{R}} & (\mathrm{i} \partial_t \psi) {u_{\varepsilon}}e^{\mathrm{i} e_{\varepsilon} t + \mathrm{i} \beta S[|\varphi u_{\varepsilon} |^2]} \mathrm{d} y \\
    & =  \int_{\mathbb{R}} \mathrm{i} \left[(\partial_t \varphi){u_{\varepsilon}}e^{-\mathrm{i} e_{\varepsilon} t - \mathrm{i} \beta S[|\varphi u_{\varepsilon} |^2]}  - \mathrm{i} \left( e_{\varepsilon}  + \beta \partial_t S[|\varphi u_{\varepsilon} |^2]\right) \varphi {u_{\varepsilon}}e^{-\mathrm{i} e_{\varepsilon} t - \mathrm{i} \beta S[|\varphi u_{\varepsilon} |^2]} \right]{u_{\varepsilon}}e^{\mathrm{i} e_{\varepsilon} t + \mathrm{i} \beta S[|\varphi u_{\varepsilon} |^2]} \mathrm{d} y \\
    & = \mathrm{i} \partial_t \varphi + e_{\varepsilon} \varphi + \beta  \varphi \int_{\mathbb{R} } \partial_t (S[|\varphi u_{\varepsilon} |^2] ) u_{\varepsilon}^2 \mathrm{d} y \\
    & = \mathrm{i} \partial_t \varphi + e_{\varepsilon} \varphi + \beta  \varphi \int_{\mathbb{R} } \partial_t (S * (|\varphi u_{\varepsilon} |^2) ) u_{\varepsilon}^2 \mathrm{d} y \\
    & = \mathrm{i} \partial_t \varphi + e_{\varepsilon} \varphi + \beta  \varphi \int_{\mathbb{R} } S * (\partial_t (|\varphi u_{\varepsilon} |^2) ) u_{\varepsilon}^2 \mathrm{d} y.
\end{split}
\end{equation}
Then using the continuity equation \eqref{cteq} for $\psi = |\varphi u_{\varepsilon} |^2 $, we obtain
\begin{equation}\nonumber
    \begin{split}
    \int_{\mathbb{R}}  (\mathrm{i} \partial_t \psi) {u_{\varepsilon}}e^{\mathrm{i} e_{\varepsilon} t + \mathrm{i} \beta S[|\varphi u_{\varepsilon} |^2]} \mathrm{d} y & = \mathrm{i} \partial_t \varphi + e_{\varepsilon} \varphi + \beta  \varphi \int_{\mathbb{R} } S * (-2 \nabla_{\mathbf{x}}\cdot \mathbf{J}_{\beta  \mathbf{A}[|\psi|^2]}(\psi) ) u_{\varepsilon}^2 \mathrm{d} y \\
    & = \mathrm{i} \partial_t \varphi + e_{\varepsilon} \varphi - 2\beta  \varphi \int_{\mathbb{R} }\left( (\nabla_{\mathbf{x}}S) * \mathbf{J}_{\beta  \mathbf{A}[|\psi|^2]}(\psi) \right) u_{\varepsilon}^2 \mathrm{d} y.
    \end{split}
\end{equation}

\noindent$\bullet$ {\bf{The current part.}} Integrating out $y$, we find
\begin{multline}\nonumber
        \int_{\mathbb{R}} \left[  - 2\beta (\nabla_{\mathbf{x}}^{\perp} \omega_0) * \mathbf{J}_{\beta \mathbf{A} [|\psi|^2]}(\psi) \psi \right] {u_{\varepsilon}}e^{\mathrm{i} e_{\varepsilon} t + \mathrm{i} \beta S[|\varphi u_{\varepsilon} |^2]} \mathrm{d} y   = - 2\beta  \varphi \int_{\mathbb{R}}   \left( (\nabla_{\mathbf{x}}^{\perp} \omega_0) * \mathbf{J}_{\beta \mathbf{A} [|\psi|^2]}(\psi) \right)   u_{\varepsilon}^2 \mathrm{d} y  \\
         = - 2\beta  \varphi \int_{\mathbb{R}}   \left( (\nabla_{\mathbf{x}}S +\mathbf{T}_0) * \mathbf{J}_{\beta \mathbf{A} [|\psi|^2]}(\psi) \right)   u_{\varepsilon}^2 \mathrm{d} y,
\end{multline}
which cancels with a term found above in the time derivative part. 
Only the term convoluted with $\mathbf{T}_0$ remains.\\

\noindent$\bullet$ {\bf{The potential part.}} With direct calculations, we obtain the following:
\begin{equation}\nonumber
\begin{split}
   \int_{\mathbb{R}}( V_{\varepsilon} \psi) {u_{\varepsilon}}e^{\mathrm{i} e_{\varepsilon} t + \mathrm{i} \beta S[|\varphi u_{\varepsilon} |^2]} \mathrm{d} y & = \int_{\mathbb{R}}\left( |x|^2 + \frac{1}{\varepsilon^2} |y|^2 \right) \varphi u_{\varepsilon}  e^{-\mathrm{i} e_{\varepsilon} t - \mathrm{i} \beta S[|\varphi u_{\varepsilon} |^2]} {u_{\varepsilon}} e^{\mathrm{i} e_{\varepsilon} t + \mathrm{i} \beta S[|\varphi u_{\varepsilon} |^2]} \mathrm{d} y\\
   & = |x|^2 \varphi + \varphi \int_{\mathbb{R}} \frac{1}{\varepsilon^2}|y|^2 u_{\varepsilon}^2 \mathrm{d} y.
\end{split}
\end{equation}

\noindent$\bullet$ {{\bf{The kinetic part.}} Thanks to a change of gauge, we have
\begin{equation}\nonumber
    \begin{split}
        \int_{\mathbb{R} } & \left[ \left( -\mathrm{i}\nabla_{\mathbf{x}}+\beta \mathbf{A}[|\psi|^2] \right)^2 \psi \right] {u_{\varepsilon}}e^{\mathrm{i} e_{\varepsilon} t + \mathrm{i} \beta S[|\varphi u_{\varepsilon} |^2]} \mathrm{d} y \\
        & \qquad \qquad =  \int_{\mathbb{R} } \left[  e^{- \mathrm{i} e_{\varepsilon} t -\mathrm{i}\beta S[|\varphi u_{\varepsilon}|^2]} \left( - \mathrm{i} \nabla_{\mathbf{x}} +  \beta \mathbf{T}[|\varphi u_{\varepsilon}|^2] \right)^2 (\varphi u_{\varepsilon })\right] {u_{\varepsilon}}e^{\mathrm{i} e_{\varepsilon} t + \mathrm{i} \beta S[|\varphi u_{\varepsilon} |^2]} \mathrm{d} y   \\
        & \qquad \qquad  = \int_{\mathbb{R} } \left[ \left( - \mathrm{i} \nabla_{\mathbf{x}} +  \beta \mathbf{T}[|\varphi u_{\varepsilon}|^2] \right)^2 (\varphi u_{\varepsilon })\right] {u_{\varepsilon}} \mathrm{d} y.
    \end{split}
\end{equation}
Then using the properties of $f$ mentioned in \eqref{fProperties}, we obtain 
\begin{equation}\nonumber
\begin{split}
        \int_{\mathbb{R} } & \left[ \left( -\mathrm{i}\nabla_{\mathbf{x}}+\beta \mathbf{A}[|\psi|^2] \right)^2 \psi \right]  {u_{\varepsilon}}e^{\mathrm{i} e_{\varepsilon} t + \mathrm{i} \beta S[|\varphi u_{\varepsilon} |^2]} \mathrm{d} y \\
        & = \int_{\mathbb{R} } \left[ \left( \mathrm{i}\partial_x + \pi \beta f(y) |\varphi(t,x)|^2 \right)^2 (\varphi(t,x) u_{\varepsilon }(y))  - \partial^2_y  (\varphi(t,x) u_{\varepsilon }(y)) \right] {u_{\varepsilon}(y)} \mathrm{d} y \\
        & = \int_{\mathbb{R} } \left[ (-\partial_x^2 \varphi) u_{\varepsilon 
        } + \mathrm{i} \pi \beta f |\varphi|^2 (\partial_x \varphi) u_{\varepsilon} + \mathrm{i} \pi \beta f \partial_x (|\varphi|^2 \varphi) u_{\varepsilon} +  \pi^2 \beta^2 f^2 |\varphi|^4 \varphi u_{\varepsilon} + \varphi ( - \partial^2_y   u_{\varepsilon }) \right] {u_{\varepsilon}} \mathrm{d} y  \\
        & = - \partial_x^2 \varphi + \frac{1}{3}\pi^2 \beta^2 |\varphi|^4 \varphi  + \varphi \int_{\mathbb{R} } ( - \partial^2_y   u_{\varepsilon }) u_{\varepsilon} \mathrm{d} y.
\end{split}
\end{equation}

Adding all the parts above, and using 
\begin{equation}\nonumber
    - \partial_y^2 u_{\varepsilon} + \frac{1}{\varepsilon^2}|y|^2 u_{\varepsilon} = e_{\varepsilon} u_{\varepsilon},
\end{equation}
it turns out that 
\begin{equation}\nonumber
    \mathrm{i} \partial_t \varphi =  - \partial_x^2 \varphi + \frac{1}{3}\pi^2 \beta^2 |\varphi|^4 \varphi +|x|^2 \varphi - 2\beta  \varphi \int_{\mathbb{R}}   \left( \mathbf{T}_0 * \mathbf{J}_{\beta \mathbf{A} [|\psi|^2]}(\psi) \right)   u_{\varepsilon}^2 \mathrm{d} y .
\end{equation}
To complete the calculation, we have the following claim.\\

{\noindent \bf Claim:} $$ - 2\beta  \varphi \int_{\mathbb{R}}   \left( \mathbf{T}_0 * \mathbf{J}_{\beta \mathbf{A} [|\psi|^2]}(\psi)  \right)   u_{\varepsilon}^2 \mathrm{d} y  =  \frac{2}{3}\pi^2 \beta^2 |\varphi|^4 \varphi. $$

\begin{proof}
Using a change of gauge, we obtain 
\begin{equation} \nonumber
    \mathbf{J}_{\beta \mathbf{A} [|\psi|^2]}(\psi)   = \mathbf{J}_{\beta \mathbf{T}[|\varphi u_{\varepsilon}|^2]}(\varphi u_{\varepsilon}).
\end{equation}
With the same definition as in \eqref{def-f}, we have
\begin{equation}\nonumber
   \mathbf{T}[|\varphi u_{\varepsilon}|^2] (t,x,y) =  \left( \begin{matrix}
       -\pi f(y) |\varphi(t,x)|^2 \\
       0
   \end{matrix} \right).
\end{equation}
Direct calculations for the current give \footnote{$c.c.$ stands for ``complex conjugate''.}
\begin{equation}\nonumber
    \begin{split}
        \mathbf{J}_{\beta \mathbf{T}[|\varphi u_{\varepsilon}|^2]}(\varphi u_{\varepsilon}) & = \frac{1}{2} \left[ \overline{(\varphi u_{\varepsilon})} \left( -\mathrm{i} \nabla_{\mathbf{x}} + \beta \mathbf{T}[|\varphi u_{\varepsilon}|^2] \right) (\varphi u_{\varepsilon})  + c. c.\right] \\
        & = \frac{1}{2} \left[ (\overline{\varphi} u_{\varepsilon}) \left( \begin{matrix}
             -\mathrm{i} \partial_x  - \pi \beta  f |\varphi|^2  \\
              -\mathrm{i} \partial_y
        \end{matrix}\right) (\varphi u_{\varepsilon})  + c. c. \right] \\
        & = \frac{1}{2} \left[ (\overline{\varphi} u_{\varepsilon}) \left( \begin{matrix}
             -\mathrm{i} (\partial_x \varphi) u_{\varepsilon}  - \pi \beta  f |\varphi|^2 \varphi u_{\varepsilon}  \\
              -\mathrm{i} \varphi \partial_y  u_{\varepsilon}
        \end{matrix}\right)  +  c. c.  \right] \\
        & = \frac{1}{2} \left( \begin{matrix}
            - \mathrm{i} \overline{\varphi} \partial_x \varphi +\mathrm{i} \varphi \partial_x \overline{\varphi} - 2 \pi \beta  f |\varphi|^4 \\
            0
        \end{matrix} \right) u_{\varepsilon}^2.
    \end{split}
\end{equation}
Then we have
\begin{equation} \nonumber
    \begin{split}
        - 2\beta  \varphi \int_{\mathbb{R}}  &   \left( \mathbf{T}_0 * \mathbf{J}_{\beta \mathbf{T}[|\varphi u_{\varepsilon}|^2]}(\varphi u_{\varepsilon}) \right)   u_{\varepsilon}^2 \mathrm{d} y \\
        & = \pi \beta \varphi \int_{\mathbb{R}}\left(  \int_{\mathbb{R}} \operatorname{sgn}(y-y') \left( - \mathrm{i} \overline{\varphi} \partial_x \varphi +\mathrm{i} \varphi \partial_x \overline{\varphi} - 2 \pi \beta  f(y') |\varphi|^4 \right) u_{\varepsilon}^2(y') \mathrm{d} y'  \right) u_{\varepsilon}^2 (y)  \mathrm{d} y \\
        & = \pi \beta \varphi \int_{\mathbb{R}}\left(  \int_{\mathbb{R}} \operatorname{sgn}(y'-y) \left( \mathrm{i} \overline{\varphi} \partial_x \varphi  -  \mathrm{i} \varphi \partial_x \overline{\varphi} + 2 \pi \beta  f(y') |\varphi|^4 \right) u_{\varepsilon}^2(y) \mathrm{d} y  \right) u_{\varepsilon}^2 (y')  \mathrm{d} y'  \\
        & = \pi \beta \varphi \int_{\mathbb{R}}\left( \mathrm{i} \overline{\varphi} \partial_x \varphi  -  \mathrm{i} \varphi \partial_x \overline{\varphi} + 2 \pi \beta  f(y') |\varphi|^4 \right) f(y') u_{\varepsilon}^2 (y')  \mathrm{d} y' \\
        & = \frac{2}{3}\pi^2 \beta^2 |\varphi|^4 \varphi, 
    \end{split}
\end{equation}
where the last step follows from the properties of $f$ mentioned in \eqref{fProperties}.
\end{proof}

This concludes the calculation, showing that if the 2D time evolution solution decouples, the part in the loose direction should be described by the quintic NLS equation. Therein, the properties of $f$ mentioned in \eqref{fProperties} play an important role.\\

\subsection{Dimensional reduction} \label{dynamicsproof}

In this section, we prove Theorem \ref{dynamicresult}. \\

We will use the orthogonal projection on $u_1$, the ground state of $H_y$ as in \eqref{def-u}:

\begin{definition}[Projection on $u_1$]
For a function $ \Psi: [0,T_0] \times\mathbb{R}^2 \to \mathbb{C} $, we define its projection on $u_1$ by
\begin{equation} \nonumber
    (\Pi_1 \Psi) (t,x,y) \colonequals \left( \int_{\mathbb{R}} \Psi(t,x,y') u_1(y') \mathrm{d} y' \right) u_1(y).
\end{equation}\end{definition}

The projection of the solution $\phi_{\varepsilon}$ to the IVP \eqref{rescaledPDE} can be rewritten as
\begin{equation}\nonumber
    (\Pi_1 \phi_{\varepsilon} )(t,x,y) = e^{-\mathrm{i} \frac{t}{\varepsilon}} \varphi_{\varepsilon} (t,x) u_1(y)
\end{equation}
with
\begin{equation} \label{def-vphi-e}
    \varphi_{\varepsilon}(t,x) \colonequals e^{\mathrm{i}\frac{t}{\varepsilon}} \int_{\mathbb{R}} \phi_{\varepsilon}(t,x,y) u_1(y) \mathrm{d} y,
\end{equation}
then the 2D dynamics \eqref{rescaledPDE} becomes
\begin{equation}\label{1DPDEwithepsilon}
    \begin{cases}
        \mathrm{i} \partial_t \varphi_{\varepsilon} = H_x \varphi_{\varepsilon} + F\left( \frac{t}{\varepsilon}, \phi_{\varepsilon} \right) \\
        \varphi_{\varepsilon}(0,x) = \varphi_0 (x)
    \end{cases}
\end{equation}
with
\begin{equation}\nonumber
    F (\theta, \Psi) = e^{\mathrm{i} \theta } \int_{\mathbb{R}}f[\Psi](t,x,y) u_1(y) \mathrm{d} y,
\end{equation}
where $f[\cdot]$ is as in \eqref{f[]}.
We now focus on the 1D dynamics \eqref{1DPDEwithepsilon} and prove that as $\varepsilon \to 0$ its solution $\varphi_{\varepsilon}$ converges to the solution $\varphi$ to the IVP \eqref{1Divp}, which is the quintic NLS equation. 

\begin{theorem}[Dynamics restatement]\label{difference}
    Let $\varphi_{\varepsilon}$ and $\varphi$ be the solutions to the initial value problems \eqref{1DPDEwithepsilon} and \eqref{1Divp}, respectively. Then we have
    \begin{equation}\nonumber
    \sup_{T\in[0,T_0]}\norm{\LR{ \varphi_{\varepsilon} - \varphi}(T,\cdot)}_{L^2_x}^{}  \le C_{T_0} \varepsilon^{\frac{1}{4}}
    \end{equation}
    for a constant $C_{T_0}$ independent of $\varepsilon$.
\end{theorem}

We define the difference
\begin{equation}\nonumber
    \chi_{\varepsilon} \colonequals \varphi_{\varepsilon} - \varphi,
\end{equation}
then it satisfies the initial value problem
\begin{equation} \label{diffchi}
    \begin{cases}
        \mathrm{i} \partial_t \chi_{\varepsilon} = H_x \chi_{\varepsilon} + F\left( \frac{t}{\varepsilon}, \phi_{\varepsilon} \right) - \pi^2 \beta^2 |\varphi|^4 \varphi\\
        \chi_{\varepsilon}(0,x) = 0
    \end{cases}.
\end{equation}
Rephrasing the calculations in Section \ref{trialPDE}, we have
\begin{equation}\nonumber
    \pi^2 \beta^2 |\varphi|^4 \varphi = F \left( \frac{t}{\varepsilon}, e^{-\mathrm{i} \frac{t}{\varepsilon}}\varphi u_1 \right) .
\end{equation}
Then the IVP \eqref{diffchi} can be rewritten as
\begin{equation}\nonumber
    \begin{cases}
        \mathrm{i} \partial_t \chi_{\varepsilon} = H_x \chi_{\varepsilon} + h_{\varepsilon} + g_{\varepsilon} \\
        \chi_{\varepsilon}(0,x) = 0
    \end{cases},
\end{equation}
where 
\begin{equation}\nonumber
    h_{\varepsilon} = F\left( \frac{t}{\varepsilon}, \phi_{\varepsilon} \right)  - F\left( \frac{t}{\varepsilon}, \Pi_1 \phi_{\varepsilon} \right)  = e^{\mathrm{i} \frac{t}{\varepsilon} } \int_{\mathbb{R}}\left( f[\phi_{\varepsilon}]- f[\Pi_1 \phi_{\varepsilon}]\right) (t,x,y)  u_1(y) \mathrm{d} y
\end{equation}
and
\begin{equation}\nonumber
    g_{\varepsilon}  = F\left( \frac{t}{\varepsilon}, \Pi_1 \phi_{\varepsilon} \right) - F \left( \frac{t}{\varepsilon}, e^{-\mathrm{i} \frac{t}{\varepsilon}}\varphi u_1 \right)  = \pi^2 \beta^2 |\varphi_{\varepsilon}|^4 \varphi_{\varepsilon} - \pi^2 \beta^2 |\varphi|^4 \varphi.
\end{equation}
To prove that the difference $\chi_{\varepsilon}$ vanishes as $\varepsilon \to 0 $, we need estimates on $h_{\varepsilon}$ and $g_{\varepsilon}$. 
In the sequel, the estimate on $h_{\varepsilon}$ is the more challenging one, which requires the uniform $H^2$ well-posedness stated in Assumption \ref{assumption}. 
We will use standard $L_t^qL^p_{\bf{x}}$ norms as follows.
\begin{definition}[$L_t^qL^p_{\bf{x}}$ norms]
For a function $\Psi:[0,T]\times\mathbb{R}^n \to \mathbb{C}$,  \, $T\in (0,T_0]$ and $p,q \ge 1$, we define the norm
\begin{equation}\nonumber
    \norm{\Psi}_{L^{q}_{0\le t\le T}(L^p_{\bf{x}})} = \left( \int_{0}^T\norm{\Psi(t,\cdot)}^{q}_{L^p}\mathrm{d} t \right)^{\frac{1}{q}} , \quad  \norm{\Psi}_{L^{\infty}_{0\le t\le T}(L^p_{\bf{x}})} = 
    \operatorname{esssup}_{0\le t\le T}\norm{\Psi(t,\cdot)}^{}_{L^p}.
\end{equation}
For $n=2$, we define the anisotropic norm
\begin{equation}\nonumber
    \norm{\Psi(t,\cdot)}_{L^{q}_{x}(L^p_{y})} = \left( \int_{\mathbb{R}}\norm{\Psi(t,x,\cdot)}^{q}_{L^p}\mathrm{d} x \right)^{\frac{1}{q}} \quad  \norm{\Psi(t,\cdot)}_{L^{\infty}_{x} (L^p_{y})} = 
    \operatorname{esssup}_{x\in \mathbb{R}}\norm{\Psi(t,x,\cdot)}^{}_{L^p}.
\end{equation}
\end{definition}

\begin{lemma}[$H^1(\mathbb{R}^2) \subset L^q_x(L^2_y)$ for $q \ge 2$]\label{lemma-embedding}
Let $\Psi \in H^1(\mathbb{R}^2)$ and $q\ge 2$. Then
\begin{equation}\nonumber
    \norm{\Psi}_{L^{q}_{x}(L^2_{y})} \le C_q \norm{\Psi}_{H^1}
\end{equation}
for some constant $C_q$ independent of $\Psi$.
\end{lemma}
\begin{proof}
    Define 
    \begin{equation}\nonumber
        f_{\Psi}(x) \colonequals \norm{\Psi(x,\cdot)}^{2}_{L^2},
    \end{equation}
    which lies in $L^1(\mathbb{R})$ because of $\Psi \in L^2(\mathbb{R}^2)$.
    With direct calculations, we have
    \begin{equation}\nonumber
        \abs{f'_{\Psi}(x) } = \abs{\int_{\mathbb{R}} \partial_x |\Psi(x,y)|^2\differential y } \le  \int_{\mathbb{R}} \LR{|\Psi(x,y)|^2 + \abs{\partial_x \Psi(x,y)}^2}\differential y = f_{\Psi}(x) + f_{\partial_x \Psi}(x) \in L^1(\mathbb{R}),
    \end{equation}
    which implies $f_{\Psi} \in W^{1,1}(\mathbb{R})$.
    Applying the Sobolev embedding $W^{1,1}(\mathbb{R}) \subset L^p(\mathbb{R})$ for $p\ge 1$, it turns out that
    \begin{equation}\nonumber
        \norm{\Psi}^q_{L^q_x(L^2_y)} = \norm{f_{\Psi}}^{\frac{q}{2}}_{L^{\frac{q}{2}}} \le C_q \norm{f_{\Psi}}^{\frac{q}{2}}_{W^{1,1}}  \le C_q  \norm{\Psi}_{H^1}^q,
    \end{equation}
    which completes the proof.
\end{proof}

\begin{prop}[Estimate on $g_{\varepsilon}$] \label{estimateg}
With the notation above, we have
    \begin{equation}\nonumber
        \norm{g_{\varepsilon}}_{L^{\frac{4}{3}}_{0\le t\le T}(L^1_x)} \le C \norm{\chi_{\varepsilon}}_{L^{\frac{4}{3}}_{0\le t\le T}(L^2_x)}
    \end{equation}
for some constant $C$ independent of $\varepsilon$ and $T$.
\end{prop}
\begin{proof}
Clearly, we have
    \begin{equation}\nonumber
        |g_{\varepsilon}| \le \frac{5}{2} \pi^2 \beta^2 \left( |\varphi_{\varepsilon}|^4 + |\varphi|^4 \right) |\chi_{\varepsilon}|.
    \end{equation}
With a direct calculation, we find
    \begin{equation}\nonumber
        \begin{split}
        \norm{g_{\varepsilon}}_{L^{\frac{4}{3}}_{0\le t\le T}(L^1_x)} = \left( \int_{0}^T\norm{g_{\varepsilon}(t,\cdot)}^{\frac{4}{3}}_{L^1}\mathrm{d} t \right)^{\frac{3}{4}}  & \le C \left( \int_{0}^T\norm{\left( |\varphi_{\varepsilon}(t,\cdot)|^4 + |\varphi(t,\cdot)|^4 \right) |\chi_{\varepsilon}(t,\cdot)|}^{\frac{4}{3}}_{L^1}\mathrm{d} t \right)^{\frac{3}{4}}\\
        \text{(H\"older's inequality)} \quad & \le C \left( \int_{0}^T\norm{ |\varphi_{\varepsilon}(t,\cdot)|^4 + |\varphi(t,\cdot)|^4 }^{\frac{4}{3}}_{L^2} \norm{\chi_{\varepsilon}(t,\cdot)}^{\frac{4}{3}}_{L^2}\mathrm{d} t \right)^{\frac{3}{4}} \\
        & \le C \left( \int_{0}^T \left( \norm{ \varphi_{\varepsilon}(t,\cdot)}^{\frac{16}{3}}_{L^8}  + \norm{\varphi(t,\cdot) }^{\frac{16}{3}}_{L^8} \right) \norm{\chi_{\varepsilon}(t,\cdot)}^{\frac{4}{3}}_{L^2}\mathrm{d} t \right)^{\frac{3}{4}} \text{}\\
        & \le C  \left( \norm{ \varphi_{\varepsilon}}^{4}_{L^{\infty}_{0\le t\le T}(L^8_x)}  + \norm{\varphi }^{4}_{L^{\infty}_{0\le t\le T}(L^8_x)} \right)\norm{\chi_{\varepsilon}}_{L^{\frac{4}{3}}_{0\le t\le T}(L^2_x)}.
        \end{split}
    \end{equation}
Recall the definition of $\varphi_{\varepsilon}$ in \eqref{def-vphi-e}. Applying Jensen's inequality, it follows that
\begin{equation}\nonumber
    |\varphi_{\varepsilon}(t,x)| \le 2 \norm{\phi_{\varepsilon}(t,x,\cdot)}_{L^2}.
\end{equation}
Then, using Lemma \ref{lemma-embedding} above, we obtain
\begin{equation}\nonumber
    \norm{ \varphi_{\varepsilon}(t,\cdot)}_{L^8} \le  C \norm{ |\phi_{\varepsilon}|(t,\cdot)}_{H^1} .
\end{equation}
For $\varphi$, using the Sobolev embedding $H^1(\mathbb{R}) \subset L^{\infty}(\mathbb{R})$, we have
\begin{equation}\nonumber
    \norm{\varphi(t,\cdot)}_{L^8} \le \norm{\varphi(t,\cdot)}_{L^{\infty}}^{\frac{3}{4}} \norm{\varphi(t,\cdot)}_{L^2}^{\frac{1}{4}} \le C \norm{\partial_x\varphi(t,\cdot)}_{L^2}^{\frac{3}{4}} \norm{\varphi(t,\cdot)}_{L^2}^{\frac{1}{4}}.
\end{equation}
Thanks to the conservation laws for mass and for energy, $\norm{ |\phi_{\varepsilon}|(t,\cdot)}_{H^1} $ and $\norm{\varphi(t,\cdot)}_{H^1} $ are bounded from above by a constant independent of $\varepsilon$ and $t$, which completes the proof.
\end{proof}

To estimate $h_{\varepsilon}$, we need the following two lemmas.

\begin{lemma}[Uniform bound in anisotropic spaces]\label{bound}
    Let $\phi_{\varepsilon}$ be the solution to the IVP \eqref{rescaledPDE}. Then $\phi_{\varepsilon}$ and $\Pi_1\phi_{\varepsilon}$ are both uniformly bounded with respect to $\varepsilon$ in $L^{\infty}( [0,T_0], L^p_x(L^2_y))$ for $p\ge 2$. And with Assumption \ref{assumption}, $\partial_x \phi_{\varepsilon}$ and $\partial_x \Pi_1\phi_{\varepsilon}$  are both uniformly bounded with respect to $\varepsilon \in ( 0, \varepsilon_0]$ in $\mathcal{C}( [0,T_0], L^p_x(L^2_y))$ for $p\ge 2$.
\end{lemma}
\begin{proof}
    Thanks to the conservation laws for mass and for energy, $\norm{ |\phi_{\varepsilon}|(t,\cdot)}_{H^1} $ is bounded from above by a constant independent of $\varepsilon$ and $t$, which implies that $|\phi_{\varepsilon}|$ is uniformly bounded with respect to $\varepsilon$ in $L^{\infty} ([0,T_0], H^1(\mathbb{R}^2))$.
    For the projection $\Pi_1$, it is straightforward to have the following properties:
    \begin{equation}\nonumber
        \norm{\Pi_1}_{L^2\to L^2} \le 1, \quad \partial_x \Pi_1 = \Pi_1 \partial_x.
    \end{equation}
    Hence, the $H^1$ uniform boundedness for $\phi_{\varepsilon}$ above implies the same $H^1$ uniform boundedness for $\Pi_1 \phi_{\varepsilon}$.     
    Applying Lemma \ref{lemma-embedding}, these uniform bounds are also true in the anisotropic space $L^{\infty} ( [0,T_0], L^p_x(L^2_y))$ for $p\ge 2$.
    Assumption \ref{assumption} gives the $H^2$ uniform boundedness for $\phi_{\varepsilon}$, which implies that $\norm{\partial_x\phi_{\varepsilon}(t,\cdot)}_{H^1}$ and $\norm{\partial_x\Pi_1\phi_{\varepsilon}(t,\cdot)}_{H^1}$ are uniformly bounded from above with respect to $\varepsilon\in(0,\varepsilon_0]$ and $t\in[0,T_0]$.
    Again using Lemma \ref{lemma-embedding}, we conclude the proof.
\end{proof}

\begin{lemma}[Projection estimate] \label{diffphiestim}
Let $\phi_{\varepsilon}$ be the solution to the IVP \eqref{rescaledPDE},
we have
\begin{equation}\nonumber
    \norm{\phi_{\varepsilon} - \Pi_1 \phi_{\varepsilon}}_{L^{\infty}_{0\le t\le T_0}(L^2_{\mathbf{x}})}^2 \le C {\varepsilon}
\end{equation}
for some constant $C$ independent of $\varepsilon$ and $T_0$.
\end{lemma}
\begin{proof}
Rewrite $\phi_{\varepsilon}$ as
    \begin{equation}\nonumber
        \phi_{\varepsilon} (t,x,y) = \sum_{k=1}^{\infty} \phi_{\varepsilon,k} (t,x) u_{1,k} (y),
    \end{equation}
where $u_{1,k}$ is the $k$th-eigenfunction of $H_y$ with respect to the $k$-th eigenvalue $\lambda_k$  and 
    \begin{equation}\nonumber
         \phi_{\varepsilon,k} (t,x) = \int_{\mathbb{R}} \phi_{\varepsilon}(t,x,y) u_{1,k}(y) \mathrm{d} y .
    \end{equation}
Clearly, 
\begin{equation}\nonumber
    u_{1,1} = u_1, \quad \lambda_1 = 1,\quad \lambda_{k+1} > \lambda_k >\lambda_1 \text{ for } k\ge 2, 
\end{equation}
\begin{equation}\nonumber
    \phi_{\varepsilon,1} (0,x) = \varphi_0 (x), \quad  \phi_{\varepsilon,k} (0,x) = 0 \text{ for } k \ge 2.
\end{equation}
Since the dynamics conserves the $L^2$ mass, we have
\begin{equation}\label{CSVofMass}
    1=\norm{\varphi_0}_{L^2}^2 = \norm{\phi_{\varepsilon}(0,\cdot)}_{L^2}^2 = \norm{\phi_{\varepsilon}(t,\cdot)}_{L^2}^2 = \sum_{k=1 }^{\infty}\norm{\phi_{\varepsilon,k}(t,\cdot)}_{L^2}^2 .
\end{equation}
The energy corresponding to \eqref{rescaledPDE} is
\begin{equation}\nonumber
    \tilde{\mathcal{E}}^{2\mathrm{D}}_{\varepsilon} (\phi_{\varepsilon})  = \frac{1}{\varepsilon} \langle\phi_{\varepsilon},H_y \phi_{\varepsilon}\rangle + \int_{\mathbb{R}^2} |(-\mathrm{i}\partial_{x} + \beta (\mathbf{T}_0)_x * |\phi_{\varepsilon}|^2) \phi_{\varepsilon}|^2 + \int_{\mathbb{R}^2} |x|^2|\phi_{\varepsilon}|^2.
\end{equation}
With direct calculations, we have
\begin{equation}\label{tenergy}
    \langle\phi_{\varepsilon}(t,\cdot) , H_y \phi_{\varepsilon}(t,\cdot) \rangle = \sum_{k=1}^{\infty} \lambda_k \norm{\phi_{\varepsilon,k}(t,\cdot)}_{L^2}^2 = \sum_{k=2}^{\infty} (\lambda_k - 1) \norm{\phi_{\varepsilon,k}(t,\cdot)}_{L^2}^2 + \sum_{k=1}^{\infty} \norm{\phi_{\varepsilon,k}(t,\cdot)}_{L^2}^2
\end{equation}
and 
\begin{equation}\label{0energy}
    \langle\phi_{\varepsilon}(0,\cdot) , H_y \phi_{\varepsilon}(0,\cdot) \rangle = \norm{\varphi_0}_{L^2}^2 = 1.
\end{equation}
Using the conservation law for energy and \eqref{CSVofMass}, \eqref{tenergy}, \eqref{0energy}, we obtain
\begin{multline}\nonumber
    \frac{1}{\varepsilon} + \mathcal{E}^{1\mathrm{D}}(\varphi_0) = \tilde{\mathcal{E}}^{2\mathrm{D}}_{\varepsilon} (\phi_{\varepsilon}(0,\cdot))   = \tilde{\mathcal{E}}^{2\mathrm{D}}_{\varepsilon} (\phi_{\varepsilon}(t,\cdot)) = \frac{1}{\varepsilon}  + \frac{1}{\varepsilon}  \sum_{k=2}^{\infty} (\lambda_k - 1) \norm{\phi_{\varepsilon,k}(t,\cdot)}_{L^2}^2  \\ +  \int_{\mathbb{R}^2} |(-\mathrm{i}\partial_{x} + \beta (\mathbf{T}_0)_x * |\phi_{\varepsilon}(t,\cdot)|^2) \phi_{\varepsilon}(t,\cdot)|^2 + \int_{\mathbb{R}^2} |x|^2|\phi_{\varepsilon}(t,\cdot)|^2,
\end{multline}
which results in 
\begin{equation}\nonumber
    \sum_{k=2}^{\infty} (\lambda_k - 1) \norm{\phi_{\varepsilon,k}(t,\cdot)}_{L^2}^2 \le C \varepsilon
\end{equation}
for some constant $C$ independent of $\varepsilon$. Hence, we get 
\begin{equation}\nonumber
    \norm{\phi_{\varepsilon} (t,\cdot) - \Pi_1 \phi_{\varepsilon}(t,\cdot) }_{L^2}^2 =  \sum_{k=2}^{\infty} \norm{\phi_{\varepsilon,k}(t,\cdot)}_{L^2}^2 \le \frac{1}{\lambda_2 - 1} \sum_{k=2}^{\infty} (\lambda_k - 1) \norm{\phi_{\varepsilon,k}(t,\cdot)}_{L^2}^2 \le C \varepsilon,
\end{equation}
which concludes the proof of Lemma \ref{diffphiestim}.
\end{proof}

\begin{prop}[Estimate on $h_{\varepsilon}$]\label{estimateh}
With Assumption \ref{assumption} and the notation above, we have
        \begin{equation}\nonumber
        \norm{h_{\varepsilon}}_{L^{\frac{4}{3}}_{0\le t\le T}(L^1_x)} \le C T^{\frac{3}{4}} \varepsilon^{\frac{1}{4}}
    \end{equation}
for some constant $C$ independent of $\varepsilon$ and $T$.
\end{prop}
\begin{proof}
Recalling the definition of $f[\cdot]$ in \eqref{f[]}:
\begin{align}
\norm{h_{\varepsilon}(t,\cdot) }_{L^1} & \le \int_{\mathbb{R}^2} \abs{ f[\phi_{\varepsilon}]- f[\Pi_1 \phi_{\varepsilon}]} (t,x,y)    u_1(y) \mathrm{d} x \mathrm{d} y   \nonumber \\
        &\le  \beta^2 \int_{\mathbb{R}^2} \abs{  \left( (\mathbf{T}_0)_x *|\phi_{\varepsilon}|^2 \right)^2 \phi_{\varepsilon}  - \left( (\mathbf{T}_0)_x *|\Pi_1\phi_{\varepsilon}|^2 \right)^2 \Pi_1 \phi_{\varepsilon}} u_1   \label{diff1}\\
        &\quad  + \beta \int_{\mathbb{R}^2} \abs{   \left( (\mathbf{T}_0)_x *|\phi_{\varepsilon}|^2  \right) \partial_x \phi_{\varepsilon}  -  \left( (\mathbf{T}_0)_x *|\Pi_1\phi_{\varepsilon}|^2  \right) \partial_x \Pi_1\phi_{\varepsilon}   } u_1   \label{diff2}  \\
        &\quad  + \beta \int_{\mathbb{R}^2} \abs{  \partial_x \left[ \left( (\mathbf{T}_0)_x *|\phi_{\varepsilon}|^2 \right) \phi_{\varepsilon}  - \left( (\mathbf{T}_0)_x *|\Pi_1\phi_{\varepsilon}|^2 \right) \Pi_1\phi_{\varepsilon}   \right] } u_1    \label{diff3}  \\
        &\quad  + 2 \beta \int_{\mathbb{R}^2} \abs{  \left(    (\mathbf{T}_0)_x * (\mathbf{J}_{\mathbf{0}}(\phi_{\varepsilon}) )_x  \right) \phi_{\varepsilon} -    \left(  (\mathbf{T}_0)_x *(\mathbf{J}_{\mathbf{0}}(\Pi_1\phi_{\varepsilon}) )_x \right) \Pi_1\phi_{\varepsilon}} u_1    \label{diff4}  \\
        &\quad  + 2 \beta^2 \int_{\mathbb{R}^2} \Big|   (\mathbf{T}_0)_x * \left( \left( (\mathbf{T}_0)_x *|\phi_{\varepsilon}|^2  \right) |\phi_{\varepsilon}|^2 \right) \phi_{\varepsilon} \nonumber \\
        &  \qquad  \qquad  \qquad  \qquad  -   (\mathbf{T}_0)_x * \left( \left( (\mathbf{T}_0)_x *|\Pi_1\phi_{\varepsilon}|^2  \right) |\Pi_1\phi_{\varepsilon}|^2 \right) \Pi_1\phi_{\varepsilon} \Big| u_1.  \label{diff5}
\end{align}
We now estimate the right-hand side term by term.\\

\noindent$\bullet$ For the square term \eqref{diff1}, we have
\begin{equation}\nonumber
    \begin{split}
        \int_{\mathbb{R}^2} & \abs{  \left( (\mathbf{T}_0)_x *|\phi_{\varepsilon}|^2 \right)^2 \phi_{\varepsilon}  - \left( (\mathbf{T}_0)_x *|\Pi_1\phi_{\varepsilon}|^2 \right)^2 \Pi_1 \phi_{\varepsilon}} u_1 \\
        & \le \int_{\mathbb{R}^2} \left(\abs{   (\mathbf{T}_0)_x *|\phi_{\varepsilon}|^2 }^2  | \phi_{\varepsilon}  - \Pi_1 \phi_{\varepsilon}| + \abs{  \left( (\mathbf{T}_0)_x *|\phi_{\varepsilon}|^2 \right)^2 - \left( (\mathbf{T}_0)_x *|\Pi_1\phi_{\varepsilon}|^2 \right)^2 } |\Pi_1 \phi_{\varepsilon}| \right) u_1 \\
        & \le \int_{\mathbb{R}}  \norm{\abs{   (\mathbf{T}_0)_x *|\phi_{\varepsilon}|^2 }^2}_{L^2_x}  \norm{ \phi_{\varepsilon}  - \Pi_1 \phi_{\varepsilon}}_{L^2_x}  u_1 \\
        & \quad \quad  + \int_{\mathbb{R}^2} \abs{ (\mathbf{T}_0)_x * \left( |\phi_{\varepsilon}|^2  + |\Pi_1\phi_{\varepsilon}|^2\right) } \abs{ (\mathbf{T}_0)_x * \left( |\phi_{\varepsilon}|^2  - |\Pi_1\phi_{\varepsilon}|^2\right) } |\Pi_1 \phi_{\varepsilon}|   u_1 \\
        & \le  \norm{\abs{   (\mathbf{T}_0)_x *|\phi_{\varepsilon}|^2 }^2 u_1 }_{L^2_{\mathbf{x}}}  \norm{ \phi_{\varepsilon}  - \Pi_1 \phi_{\varepsilon}}_{L^2_{\mathbf{x}}}  \\
        & \quad \quad  + \int_{\mathbb{R}^2} \abs{ (\mathbf{T}_0)_x * \left( |\phi_{\varepsilon}|^2  + |\Pi_1\phi_{\varepsilon}|^2\right) } \abs{ (\mathbf{T}_0)_x * \left( \left({|\phi_{\varepsilon}|  + |\Pi_1\phi_{\varepsilon}|} \right) \left({|\phi_{\varepsilon}|  - |\Pi_1\phi_{\varepsilon}|}\right) \right) } |\Pi_1 \phi_{\varepsilon}|   u_1.
    \end{split}
\end{equation}
Recalling the definition of $(\mathbf{T}_0)_x$ in \eqref{T0x}, we have
\begin{equation}\label{Tproperty}
    \begin{split}
        \abs{(\mathbf{T}_{0})_x * \rho } (t,x,y) \le \pi \int_{\mathbb{R}}\abs{\rho(t,x,\nu)} \mathrm{d} \nu  = \pi \norm{\rho}_{L^1_y}.
    \end{split}
\end{equation}
Then the square term \eqref{diff1} can be estimated as follows:
\begin{equation}\nonumber
    \begin{split}
        \int_{\mathbb{R}^2} & \abs{  \left( (\mathbf{T}_0)_x *|\phi_{\varepsilon}|^2 \right)^2 \phi_{\varepsilon}  - \left( (\mathbf{T}_0)_x *|\Pi_1\phi_{\varepsilon}|^2 \right)^2 \Pi_1 \phi_{\varepsilon}} u_1 \\
        & \le C \norm{ \norm{\phi_{\varepsilon} }^4_{L^2_y} u_1 }_{L^2_{\mathbf{x}}}  \norm{ \phi_{\varepsilon}  - \Pi_1 \phi_{\varepsilon}}_{L^2_{\mathbf{x}}}  \\
        & \qquad \quad  + C \int_{\mathbb{R}^2} \norm{ |\phi_{\varepsilon}|^2  + |\Pi_1\phi_{\varepsilon}|^2 }_{L^1_y} \norm{ \left({|\phi_{\varepsilon}|  + |\Pi_1\phi_{\varepsilon}|} \right) \left({|\phi_{\varepsilon}|  - |\Pi_1\phi_{\varepsilon}|}\right)  }_{L^1_y} |\Pi_1 \phi_{\varepsilon}|   u_1\\
        & \le C \norm{\phi_{\varepsilon} }^4_{L^8_x(L^2_y)}  \norm{ \phi_{\varepsilon}  - \Pi_1 \phi_{\varepsilon}}_{L^2_{\mathbf{x}}}  \\
        & \qquad \quad  + C \int_{\mathbb{R}} \left({\norm{ \phi_{\varepsilon} }_{L^2_y}^2  + \norm{\Pi_1\phi_{\varepsilon} }_{L^2_y}^2 }\right)  \left({\norm{ \phi_{\varepsilon} }_{L^2_y}  + \norm{\Pi_1\phi_{\varepsilon} }_{L^2_y} }\right) \norm{{\phi_{\varepsilon}  - \Pi_1\phi_{\varepsilon}}}_{L^2_y} \norm{\Pi_1 \phi_{\varepsilon}}_{L^2_y} \\
        & \le C \norm{\phi_{\varepsilon} }^4_{L^8_x(L^2_y)}  \norm{ \phi_{\varepsilon}  - \Pi_1 \phi_{\varepsilon}}_{L^2_{\mathbf{x}}}  \\
        & \qquad \quad  + C  \norm{\norm{ \phi_{\varepsilon} }_{L^2_y}^2  + \norm{\Pi_1\phi_{\varepsilon} }_{L^2_y}^2 }_{L_x^4}  \norm{\norm{ \phi_{\varepsilon} }_{L^2_y}  + \norm{\Pi_1\phi_{\varepsilon} }_{L^2_y} }_{L_x^8} \norm{\Pi_1 \phi_{\varepsilon}}_{L^8_x(L^2_y)} \norm{{\phi_{\varepsilon}  - \Pi_1\phi_{\varepsilon}}}_{L^2_{\mathbf{x}}}  \\
        & \le C \cdot \left(  \norm{\phi_{\varepsilon} }^4_{L^8_x(L^2_y)}  + \norm{\Pi_1 \phi_{\varepsilon}}_{L^8_x(L^2_y)}^4 \right) \cdot \norm{ \phi_{\varepsilon}  - \Pi_1 \phi_{\varepsilon}}_{L^2_{\mathbf{x}}} .
    \end{split}
\end{equation}

\noindent$\bullet$ Using the property in \eqref{Tproperty}, the first cross term \eqref{diff2} becomes
\begin{equation}\nonumber
    \begin{split}
        & \int_{\mathbb{R}^2} \abs{   \left( (\mathbf{T}_0)_x *|\phi_{\varepsilon}|^2  \right) \partial_x \phi_{\varepsilon}  -  \left( (\mathbf{T}_0)_x *|\Pi_1\phi_{\varepsilon}|^2  \right) \partial_x \Pi_1\phi_{\varepsilon}   } u_1 \\
        & \le \int_{\mathbb{R}^2} \left( {\abs{   \left( (\mathbf{T}_0)_x *(|\phi_{\varepsilon}|^2  - |\Pi_1\phi_{\varepsilon}|^2)  \right) \partial_x  \phi_{\varepsilon} } + \abs{\left( (\mathbf{T}_0)_x *|\Pi_1\phi_{\varepsilon}|^2  \right) \partial_x (\phi_{\varepsilon}  - \Pi_1\phi_{\varepsilon} )  }} \right) u_1 \\
        & \le C \int_{\mathbb{R}^2} \norm{|\phi_{\varepsilon}|^2  - |\Pi_1\phi_{\varepsilon}|^2}_{L^1_y} \abs{ \partial_x \phi_{\varepsilon}   } u_1 + C \int_{\mathbb{R}^2}  \norm{    |\Pi_1\phi_{\varepsilon} |^2}_{L_y^1} \abs{\partial_x ( \phi_{\varepsilon}  -  \Pi_1 \phi_{\varepsilon})} u_1  \\
        & \le C \int_{\mathbb{R}} \left(\norm{\phi_{\varepsilon}}_{L^2_y}  + \norm{\Pi_1\phi_{\varepsilon}}_{L^2_y}\right) \norm{ \phi_{\varepsilon} -  \Pi_1\phi_{\varepsilon}}_{L^2_y} \norm{ \partial_x \phi_{\varepsilon}   }_{L^2_y} + C \int_{\mathbb{R}}  \norm{    \Pi_1\phi_{\varepsilon} }_{L_y^2}^2 \norm{\partial_x ( \phi_{\varepsilon}  -  \Pi_1 \phi_{\varepsilon})}_{L^2_y}\\
        & \le C \LR{\norm{\phi_{\varepsilon}}_{L^4_x(L^2_y)}  + \norm{\Pi_1\phi_{\varepsilon}}_{L^4_x(L^2_y)}} \norm{ \phi_{\varepsilon} -  \Pi_1\phi_{\varepsilon}}_{L^2_{\mathbf{x}}} \norm{ \partial_x \phi_{\varepsilon}   }_{L^4_x(L^2_y)} + C  {\norm{    \Pi_1\phi_{\varepsilon} }_{L^4_x(L_y^2)}^2} \norm{\partial_x ( \phi_{\varepsilon}  -  \Pi_1 \phi_{\varepsilon})}_{L^2_\mathbf{x}}.
    \end{split}
\end{equation}
With direct calculations, we find
\begin{equation}\label{dervativeestimate}
    \begin{split}
        \norm{\partial_x ( \phi_{\varepsilon}  -  \Pi_1 \phi_{\varepsilon})}_{L^2_{\mathbf{x}}}^2  & = \int_{\mathbb{R}^2} {{\partial_x ( \overline{\phi_{\varepsilon}}  -  \overline{\Pi_1 \phi_{\varepsilon}})}} {\partial_x ( \phi_{\varepsilon}  -  \Pi_1 \phi_{\varepsilon})}\\
        & = - \int_{\mathbb{R}^2} {{\partial_x^2 ( \overline{\phi_{\varepsilon}}  -  \overline{\Pi_1 \phi_{\varepsilon}})}} {( \phi_{\varepsilon}  -  \Pi_1 \phi_{\varepsilon})} \\
        & \le  \LR{\norm{\partial_x^2 {\phi_{\varepsilon}}}_{L_{\mathbf{x}}^2}  + \norm{\partial_x^2   {\Pi_1 \phi_{\varepsilon}}}_{L_{\mathbf{x}}^2} } \norm{ \phi_{\varepsilon}  -  \Pi_1 \phi_{\varepsilon}}_{L^2_{\mathbf{x}}}\\
        & \le  C {\norm{\partial_x^2 {\phi_{\varepsilon}}}_{L_{\mathbf{x}}^2} } \norm{ \phi_{\varepsilon}  -  \Pi_1 \phi_{\varepsilon}}_{L^2_{\mathbf{x}}}.
    \end{split}
\end{equation}
Then the first cross term \eqref{diff2} can be estimated as follows:
\begin{multline}\nonumber
        \int_{\mathbb{R}^2} \abs{   \left( (\mathbf{T}_0)_x *|\phi_{\varepsilon}|^2  \right) \partial_x \phi_{\varepsilon}  -  \left( (\mathbf{T}_0)_x *|\Pi_1\phi_{\varepsilon}|^2  \right) \partial_x \Pi_1\phi_{\varepsilon}   } u_1 \\
        \le C \cdot \LR{\norm{\phi_{\varepsilon}}_{L^4_x(L^2_y)}^2  + \norm{\Pi_1\phi_{\varepsilon}}^2_{L^4_x(L^2_y)} + \norm{ \partial_x \phi_{\varepsilon}   }^2_{L^4_x(L^2_y)} } \cdot  \norm{ \phi_{\varepsilon} -  \Pi_1\phi_{\varepsilon}}_{L^2_{\mathbf{x}}}  \\
        + C  \cdot  \norm{  \Pi_1\phi_{\varepsilon} }_{L^4_x(L_y^2)}^2  \cdot \norm{\partial_x^2 {\phi_{\varepsilon}}}_{L_{\mathbf{x}}^2}^{\frac{1}{2}} \cdot  \norm{ \phi_{\varepsilon}  -  \Pi_1 \phi_{\varepsilon}}_{L^2_{\mathbf{x}}}^{\frac{1}{2}}.
\end{multline}

\noindent$\bullet$ For the second cross term \eqref{diff3}, it suffices to estimate
\begin{equation}\label{diff3-}
    \begin{split}
        \int_{\mathbb{R}^2} \abs{    \left( (\mathbf{T}_0)_x * \partial_x(|\phi_{\varepsilon}|^2 )\right) \phi_{\varepsilon}  - \left( (\mathbf{T}_0)_x *\partial_x (|\Pi_1\phi_{\varepsilon}|^2) \right)\Pi_1\phi_{\varepsilon}     } u_1.
    \end{split}
\end{equation}
Applying \eqref{Tproperty} to \eqref{diff3-}, it becomes
\begin{equation}\nonumber
    \begin{split}
        \int_{\mathbb{R}^2} & \abs{    \left( (\mathbf{T}_0)_x * \partial_x(|\phi_{\varepsilon}|^2 )\right) \phi_{\varepsilon}  - \left( (\mathbf{T}_0)_x *\partial_x (|\Pi_1\phi_{\varepsilon}|^2) \right)\Pi_1\phi_{\varepsilon}     } u_1 \\
        & \le \int_{\mathbb{R}^2} \LR{\abs{ \left( (\mathbf{T}_0)_x * \partial_x(|\phi_{\varepsilon}|^2 )\right) ( \phi_{\varepsilon} - \Pi_1  \phi_{\varepsilon}) }  + \abs{\left( (\mathbf{T}_0)_x *  \partial_x (|\phi_{\varepsilon}|^2 - |\Pi_1\phi_{\varepsilon}|^2) \right) \Pi_1\phi_{\varepsilon}    }} u_1 \\
        & \le C \int_{\mathbb{R}} \norm{ \partial_x(|\phi_{\varepsilon}|^2 )}_{L^1_y} \norm{  \phi_{\varepsilon} - \Pi_1  \phi_{\varepsilon} }_{L^2_y} + C \int_{\mathbb{R}} \norm{ \partial_x (|\phi_{\varepsilon}|^2 - |\Pi_1\phi_{\varepsilon}|^2) }_{L^1_y} \norm{\Pi_1\phi_{\varepsilon} }_{L^2_y} \\
        & \le C \norm{ \partial_x(|\phi_{\varepsilon}|^2 )}_{L^2_x(L^1_y)} \norm{  \phi_{\varepsilon} - \Pi_1  \phi_{\varepsilon} }_{L^2_{\mathbf{x}}} + C\int_{\mathbb{R}} \norm{ \partial_x \LR{\phi_{\varepsilon}\overline{\phi_{\varepsilon}} - \Pi_1\phi_{\varepsilon}\overline{\Pi_1\phi_{\varepsilon}}} }_{L^1_{{y}}} \norm{\Pi_1\phi_{\varepsilon} }_{L^2_y} .
    \end{split}
\end{equation}
Since 
\begin{equation}\nonumber
    \abs{ \partial_x(|\phi_{\varepsilon}|^2 )} \le 2 \abs{\phi_{\varepsilon}}\abs{\partial_x\phi_{\varepsilon}} \le \abs{\phi_{\varepsilon}}^2 + \abs{\partial_x\phi_{\varepsilon}}^2,
\end{equation}
we have
\begin{equation}\nonumber
    \begin{split}
        \norm{ \partial_x(|\phi_{\varepsilon}|^2 )}_{L^2_x(L^1_y)} & \le \norm{\abs{\phi_{\varepsilon}}^2 + \abs{\partial_x\phi_{\varepsilon}}^2}_{L^2_x(L^1_y)} \le \norm{{\phi_{\varepsilon}}}_{L^4_x(L^2_y)}^2 + \norm{\partial_x\phi_{\varepsilon}}_{L^4_x(L^2_y)}^2.
    \end{split}
\end{equation}
On the other hand,
\begin{equation}\nonumber
    \begin{split}
        & \abs{ \partial_x \LR{\phi_{\varepsilon}\overline{\phi_{\varepsilon}} - \Pi_1\phi_{\varepsilon}\overline{\Pi_1\phi_{\varepsilon}}} } \\
        & \le  \abs{\partial_x\LR{(\phi_{\varepsilon} - \Pi_1 \phi_{\varepsilon})\overline{\phi_{\varepsilon}}}} + \abs{\partial_x\LR{{\Pi_1\phi_{\varepsilon}}(\overline{\phi_{\varepsilon}} - \overline{\Pi_1\phi_{\varepsilon}})}} \\
        & \le \abs{\overline{\phi_{\varepsilon}}\partial_x{(\phi_{\varepsilon} - \Pi_1 \phi_{\varepsilon})}} + \abs{{(\phi_{\varepsilon} - \Pi_1 \phi_{\varepsilon})\partial_x\overline{\phi_{\varepsilon}}}} + \abs{(\overline{\phi_{\varepsilon}} - \overline{\Pi_1\phi_{\varepsilon}})\partial_x{\Pi_1\phi_{\varepsilon}}}  + \abs{{{\Pi_1\phi_{\varepsilon}}\partial_x (\overline{\phi_{\varepsilon}} - \overline{\Pi_1\phi_{\varepsilon}})}},
    \end{split}
\end{equation}
hence, we have
\begin{equation}\nonumber
    \begin{split}
        \int_{\mathbb{R}} & \norm{ \partial_x \LR{\phi_{\varepsilon}\overline{\phi_{\varepsilon}} - \Pi_1\phi_{\varepsilon}\overline{\Pi_1\phi_{\varepsilon}}} }_{L^1_{{y}}} \norm{\Pi_1\phi_{\varepsilon} }_{L^2_y} \\
        & \le  \int_{\mathbb{R}}  \norm{\Pi_1\phi_{\varepsilon} }_{L^2_y}  \LR{\norm{\partial_x \phi_{\varepsilon}}_{L^2_{y}} + \norm{\partial_x \Pi_1\phi_{\varepsilon}}_{L^2_{y}} }\norm{\phi_{\varepsilon} - \Pi_1\phi_{\varepsilon} }_{L^2_{y}} \\
        & \qquad \qquad   \qquad  \qquad  \qquad +  \int_{\mathbb{R}} \norm{\Pi_1\phi_{\varepsilon} }_{L^2_y} \LR{\norm{ \phi_{\varepsilon}}_{L^2_{y}} + \norm{\Pi_1\phi_{\varepsilon}}_{L^2_{y}}}\norm{\partial_x (\phi_{\varepsilon} - \Pi_1\phi_{\varepsilon})}_{L^2_{y}}  \\
        & \le C \norm{\Pi_1\phi_{\varepsilon} }_{L^4_x(L^2_y)}\LR{\norm{\partial_x \phi_{\varepsilon}}_{L^4_x(L^2_y)}   + \norm{\partial_x \Pi_1 \phi_{\varepsilon}}_{L^4_x(L^2_y)}  }\norm{\phi_{\varepsilon} - \Pi_1\phi_{\varepsilon} }_{L^2_{\mathbf{x}}}  \\
        & \qquad \qquad  \qquad  \qquad   + C \norm{\Pi_1\phi_{\varepsilon} }_{L^4_x(L^2_y)}\LR{\norm{\phi_{\varepsilon}}_{L^4_x(L^2_y)}   + \norm{ \Pi_1 \phi_{\varepsilon}}_{L^4_x(L^2_y)}  }\norm{\partial_x (\phi_{\varepsilon} - \Pi_1\phi_{\varepsilon})}_{L^2_{\mathbf{x}}}  \\
        & \le C \LR{\norm{\Pi_1\phi_{\varepsilon} }_{L^4_x(L^2_y)}^2 + \norm{\partial_x \phi_{\varepsilon}}_{L^4_x(L^2_y)}^2   + \norm{\partial_x \Pi_1 \phi_{\varepsilon}}_{L^4_x(L^2_y)}^2  }\norm{\phi_{\varepsilon} - \Pi_1\phi_{\varepsilon} }_{L^2_{\mathbf{x}}}  \\
        & \qquad  \qquad  \qquad  + C \LR{ \norm{\phi_{\varepsilon}}_{L^4_x(L^2_y)}^2  + \norm{ \Pi_1 \phi_{\varepsilon}}_{L^4_x(L^2_y)}^2  } \norm{\partial_x^2 {\phi_{\varepsilon}}}_{L_{\mathbf{x}}^2}^{\frac{1}{2}} \norm{ \phi_{\varepsilon}  -  \Pi_1 \phi_{\varepsilon}}_{L^2_{\mathbf{x}}}^{\frac{1}{2}}.  \text{(using \eqref{dervativeestimate})}
    \end{split}
\end{equation}
Therefore, \eqref{diff3-} can be estimated as follows:
\begin{equation}\nonumber
    \begin{split}
        \int_{\mathbb{R}^2} & \abs{    \left( (\mathbf{T}_0)_x * \partial_x(|\phi_{\varepsilon}|^2 )\right) \phi_{\varepsilon}  - \left( (\mathbf{T}_0)_x *\partial_x (|\Pi_1\phi_{\varepsilon}|^2) \right)\Pi_1\phi_{\varepsilon}     } u_1 \\
        & \le C  \cdot \LR{\norm{{\phi_{\varepsilon}}}_{L^4_x(L^2_y)}^2 + \norm{\partial_x\phi_{\varepsilon}}_{L^4_x(L^2_y)}^2 + \norm{\Pi_1\phi_{\varepsilon} }_{L^4_x(L^2_y)}^2   + \norm{\partial_x \Pi_1 \phi_{\varepsilon}}_{L^4_x(L^2_y)}^2  } \cdot  \norm{  \phi_{\varepsilon} - \Pi_1  \phi_{\varepsilon} }_{L^2_{\mathbf{x}}} \\
        & \qquad\qquad \qquad\qquad \qquad + C  \cdot  \LR{ \norm{\phi_{\varepsilon}}_{L^4_x(L^2_y)}^2  + \norm{ \Pi_1 \phi_{\varepsilon}}_{L^4_x(L^2_y)}^2  }  \cdot \norm{\partial_x^2 {\phi_{\varepsilon}}}_{L_{\mathbf{x}}^2}^{\frac{1}{2}}  \cdot  \norm{ \phi_{\varepsilon}  -  \Pi_1 \phi_{\varepsilon}}_{L^2_{\mathbf{x}}}^{\frac{1}{2}}.
    \end{split}
\end{equation}

\noindent$\bullet$ With the help of \eqref{Tproperty}, the first current term \eqref{diff4} becomes
\begin{equation}\nonumber
    \begin{split}
        \int_{\mathbb{R}^2} & \abs{  \left(    (\mathbf{T}_0)_x * (\mathbf{J}_{\mathbf{0}}(\phi_{\varepsilon}) )_x  \right) \phi_{\varepsilon} -    \left(  (\mathbf{T}_0)_x *(\mathbf{J}_{\mathbf{0}}(\Pi_1\phi_{\varepsilon}) )_x \right) \Pi_1\phi_{\varepsilon}} u_1 \\
        & \le \int_{\mathbb{R}^2} \LR{\abs{  \left(    (\mathbf{T}_0)_x * (\mathbf{J}_{\mathbf{0}}(\phi_{\varepsilon}) )_x  \right) \LR{\phi_{\varepsilon} -  \Pi_1\phi_{\varepsilon}} } +  \abs{ \left(    (\mathbf{T}_0)_x * {(\mathbf{J}_{\mathbf{0}}(\phi_{\varepsilon})  - \mathbf{J}_{\mathbf{0}}(\Pi_1\phi_{\varepsilon}) )_x }\right) \Pi_1\phi_{\varepsilon}} }u_1 \\
        & \le C \int_{\mathbb{R}} \norm{    (\mathbf{J}_{\mathbf{0}}(\phi_{\varepsilon}) )_x }_{L^1_y} \norm{\phi_{\varepsilon} -  \Pi_1\phi_{\varepsilon}}_{L^2_y} +  C \int_{\mathbb{R}}  \norm{(\mathbf{J}_{\mathbf{0}}(\phi_{\varepsilon})  - \mathbf{J}_{\mathbf{0}}(\Pi_1\phi_{\varepsilon}) )_x }_{L^1_y} \norm{\Pi_1\phi_{\varepsilon}}_{L^2_y}.
    \end{split}
\end{equation}
Recalling the definition of $\mathbf{J}_0$ in \eqref{defJ}, we have
\begin{equation}\nonumber
    \begin{split}
        (\mathbf{J}_{\mathbf{0}}(\phi_{\varepsilon})  - \mathbf{J}_{\mathbf{0}}(\Pi_1\phi_{\varepsilon}) )_x & = \LR{-\frac{\mathrm{i} }{2} \overline{\phi_{\varepsilon}} \partial_x \phi_{\varepsilon} +  \frac{\mathrm{i} }{2} \overline{\Pi_1\phi_{\varepsilon}} \partial_x \Pi_1\phi_{\varepsilon}} + c.c. \\
        & = \LR{-\frac{\mathrm{i} }{2} \overline{\phi_{\varepsilon}} \partial_x \LR{\phi_{\varepsilon} - \Pi_1\phi_{\varepsilon}} - \frac{\mathrm{i} }{2} \LR{\overline{\phi_{\varepsilon}}  -\overline{\Pi_1\phi_{\varepsilon}} } \partial_x \Pi_1\phi_{\varepsilon}} + c.c.
    \end{split}
\end{equation}
Then we have
\begin{equation}\nonumber
    \begin{split}
        \norm{(\mathbf{J}_{\mathbf{0}}(\phi_{\varepsilon})  - \mathbf{J}_{\mathbf{0}}(\Pi_1\phi_{\varepsilon}) )_x }_{L^1_y} \le C \norm{\phi_{\varepsilon}}_{L^2_y} \norm{\partial_x \LR{\phi_{\varepsilon} - \Pi_1\phi_{\varepsilon}}}_{L^2_y} + C \norm{\phi_{\varepsilon} - \Pi_1\phi_{\varepsilon}}_{L^2_y} \norm{\partial_x\Pi_1\phi_{\varepsilon}}_{L^2_y}.
    \end{split}
\end{equation}
Hence, the first current term \eqref{diff4} can be estimated as follows: 
\begin{equation}\nonumber
    \begin{split}
        \int_{\mathbb{R}^2} & \abs{  \left(    (\mathbf{T}_0)_x * (\mathbf{J}_{\mathbf{0}}(\phi_{\varepsilon}) )_x  \right) \phi_{\varepsilon} -    \left(  (\mathbf{T}_0)_x *(\mathbf{J}_{\mathbf{0}}(\Pi_1\phi_{\varepsilon}) )_x \right) \Pi_1\phi_{\varepsilon}} u_1 \\
        & \le C \int_{\mathbb{R}} \norm{\phi_{\varepsilon}}_{L^2_y} \norm{\partial_x\phi_{\varepsilon}}_{L^2_y}\norm{\phi_{\varepsilon} -  \Pi_1\phi_{\varepsilon}}_{L^2_y} +  C \int_{\mathbb{R}}  \norm{\phi_{\varepsilon}}_{L^2_y} \norm{\partial_x \LR{\phi_{\varepsilon} - \Pi_1\phi_{\varepsilon}}}_{L^2_y} \norm{\Pi_1\phi_{\varepsilon}}_{L^2_y}\\
        & \qquad\qquad+ C \int_{\mathbb{R}}  \norm{\phi_{\varepsilon} - \Pi_1\phi_{\varepsilon}}_{L^2_y} \norm{\partial_x\Pi_1\phi_{\varepsilon}}_{L^2_y} \norm{\Pi_1\phi_{\varepsilon}}_{L^2_y}\\
        & \le C \LR{\norm{\phi_{\varepsilon}}_{L^4_x{(L^2_y)}} \norm{\partial_x\phi_{\varepsilon}}_{L^4_x(L^2_y)} + \norm{\Pi_1\phi_{\varepsilon}}_{L^4_x(L^2_y)} \norm{\partial_x\Pi_1\phi_{\varepsilon}}_{L^4_x(L^2_y)}} \norm{\phi_{\varepsilon} -  \Pi_1\phi_{\varepsilon}}_{L^2_{\mathbf{x}}} \\
        &\qquad\qquad +  C  \norm{\phi_{\varepsilon}}_{L^4_x(L^2_y)} \norm{\Pi_1\phi_{\varepsilon}}_{L^4_x(L^2_y)} \norm{\partial_x \LR{\phi_{\varepsilon} - \Pi_1\phi_{\varepsilon}}}_{L^2_{\mathbf{x}}}  \\
        & \le C  \cdot \LR{\norm{\phi_{\varepsilon}}_{L^4_x{(L^2_y)}}^2 + \norm{\partial_x\phi_{\varepsilon}}^2_{L^4_x(L^2_y)} + \norm{\Pi_1\phi_{\varepsilon}}^2_{L^4_x(L^2_y)} + \norm{\partial_x\Pi_1\phi_{\varepsilon}}^2_{L^4_x(L^2_y)}} \cdot  \norm{\phi_{\varepsilon} -  \Pi_1\phi_{\varepsilon}}_{L^2_{\mathbf{x}}} \\
        & \qquad\qquad +  C  \cdot \LR{\norm{\phi_{\varepsilon}}^2_{L^4_x(L^2_y)} +\norm{\Pi_1\phi_{\varepsilon}}^2_{L^4_x(L^2_y)}}  \cdot  \norm{\partial_x^2 {\phi_{\varepsilon}}}_{L_{\mathbf{x}}^2} ^{\frac{1}{2}} \cdot  \norm{ \phi_{\varepsilon}  -  \Pi_1 \phi_{\varepsilon}}_{L^2_{\mathbf{x}}}^{\frac{1}{2}}. \quad \text{(using \eqref{dervativeestimate})}
    \end{split}
\end{equation}

\noindent$\bullet$ Thanks again to \eqref{Tproperty}, the second current term \eqref{diff5} can be estimated as follows:
\begin{equation}\nonumber
    \begin{split}
        & \int_{\mathbb{R}^2} \Big|   (\mathbf{T}_0)_x * \left( \left( (\mathbf{T}_0)_x *|\phi_{\varepsilon}|^2  \right) |\phi_{\varepsilon}|^2 \right) \phi_{\varepsilon} -   (\mathbf{T}_0)_x * \left( \left( (\mathbf{T}_0)_x *|\Pi_1\phi_{\varepsilon}|^2  \right) |\Pi_1\phi_{\varepsilon}|^2 \right) \Pi_1\phi_{\varepsilon} \Big| u_1\\
        & \le  C \int_{\mathbb{R}}    \norm{\left( (\mathbf{T}_0)_x *|\phi_{\varepsilon}|^2  \right) |\phi_{\varepsilon}|^2}_{L^1_y} \norm{\phi_{\varepsilon} -  \Pi_1 \phi_{\varepsilon}}_{L^2_y} \\
        & \qquad \qquad\qquad  + C \int_{\mathbb{R}}  \norm{\left( (\mathbf{T}_0)_x *|\phi_{\varepsilon}|^2  \right) |\phi_{\varepsilon}|^2 - \left( (\mathbf{T}_0)_x *|\Pi_1\phi_{\varepsilon}|^2  \right) |\Pi_1\phi_{\varepsilon}|^2 }_{L^1_y} \norm{\Pi_1\phi_{\varepsilon} }_{L^2_y} \\
        & \le C \int_{\mathbb{R}} \norm{|\phi_{\varepsilon}|^2}_{L^1_y}^2 \norm{\phi_{\varepsilon} -  \Pi_1 \phi_{\varepsilon}}_{L^2_y}  + C \int_{\mathbb{R}}  \norm{\left( (\mathbf{T}_0)_x *|\phi_{\varepsilon}|^2  \right) \LR{|\phi_{\varepsilon}|^2 - |\Pi_1\phi_{\varepsilon}|^2} }_{L^1_y} \norm{\Pi_1\phi_{\varepsilon} }_{L^2_y} \\
        &\qquad \qquad \qquad +  C \int_{\mathbb{R}}  \norm{\LR{ (\mathbf{T}_0)_x *\LR{|\phi_{\varepsilon}|^2  - |\Pi_1\phi_{\varepsilon}|^2  }}|\Pi_1\phi_{\varepsilon}|^2 }_{L^1_y} \norm{\Pi_1\phi_{\varepsilon} }_{L^2_y} \\
        & \le  C \int_{\mathbb{R}} \norm{\phi_{\varepsilon}}_{L^2_y}^4 \norm{\phi_{\varepsilon} -  \Pi_1 \phi_{\varepsilon}}_{L^2_y}  + C \int_{\mathbb{R}} \norm{\phi_{\varepsilon}}_{L^2_y}^2  \LR{\norm{ \phi_{\varepsilon} }_{L^2_y} + \norm{\Pi_1\phi_{\varepsilon}  }_{L^2_y} }\norm{ {\phi_{\varepsilon} - \Pi_1\phi_{\varepsilon}} }_{L^2_y} \norm{\Pi_1\phi_{\varepsilon} }_{L^2_y} \\
        & \qquad \qquad\qquad +  C \int_{\mathbb{R}} \LR{\norm{ \phi_{\varepsilon} }_{L^2_y} + \norm{\Pi_1\phi_{\varepsilon}  }_{L^2_y} }\norm{ {\phi_{\varepsilon} - \Pi_1\phi_{\varepsilon}} }_{L^2_y} \norm{\Pi_1\phi_{\varepsilon} }_{L^2_y}^3 \\
        & \le C \int_{\mathbb{R}} \LR{\norm{\phi_{\varepsilon}}_{L^2_y}^4 + \norm{\Pi_1\phi_{\varepsilon} }_{L^2_y}^4} \norm{\phi_{\varepsilon} -  \Pi_1 \phi_{\varepsilon}}_{L^2_y}   \\
        & \le C \cdot \LR{\norm{\phi_{\varepsilon}}_{L^8_x(L^2_y)}^4 + \norm{\Pi_1\phi_{\varepsilon} }_{L^8_x(L^2_y)}^4}  \cdot \norm{\phi_{\varepsilon} -  \Pi_1 \phi_{\varepsilon}}_{L^2_{\mathbf{x}}} . 
    \end{split}
\end{equation}

Gathering the above estimates, we obtain
\begin{multline}\nonumber
        \norm{h_{\varepsilon}(t,\cdot)}_{L^1} \le C \cdot \Big( \norm{\phi_{\varepsilon}}_{L^8_x(L^2_y)}^4 + \norm{\Pi_1\phi_{\varepsilon} }_{L^8_x(L^2_y)}^4 +\norm{\phi_{\varepsilon}}_{L^4_x{(L^2_y)}}^2 + \norm{\partial_x\phi_{\varepsilon}}^2_{L^4_x(L^2_y)} \\
        + \norm{\Pi_1\phi_{\varepsilon}}^2_{L^4_x(L^2_y)} + \norm{\partial_x\Pi_1\phi_{\varepsilon}}^2_{L^4_x(L^2_y)} \Big) \cdot \norm{\phi_{\varepsilon} -  \Pi_1 \phi_{\varepsilon}}_{L^2_{\mathbf{x}}} \\
        +   C  \cdot \LR{\norm{\phi_{\varepsilon}}^2_{L^4_x(L^2_y)} +\norm{\Pi_1\phi_{\varepsilon}}^2_{L^4_x(L^2_y)}}  \cdot  \norm{\partial_x^2 {\phi_{\varepsilon}}}_{L_{\mathbf{x}}^2} ^{\frac{1}{2}} \cdot  \norm{ \phi_{\varepsilon}  -  \Pi_1 \phi_{\varepsilon}}_{L^2_{\mathbf{x}}}^{\frac{1}{2}}.
\end{multline}
Combining with Lemma \ref{bound} and Lemma \ref{diffphiestim}, it follows that
\begin{equation}\nonumber
    \norm{h_{\varepsilon}}_{L^{\frac{4}{3}}_{0\le t\le T}(L^1_x)} = \LR{\int_{0}^{T} \norm{h_{\varepsilon}(t,\cdot)}_{L^1}^{\frac{4}{3}} \mathrm{d} t }^{\frac{3}{4}}
    \le \LR{\int_{0}^{T} \LR{C \varepsilon^{\frac{1}{4}}}^{\frac{4}{3}} \mathrm{d} t  }^{\frac{3}{4}}
    = CT^{\frac{3}{4}} \varepsilon^{\frac{1}{4}}
\end{equation}
for some constant $C$ independent of $\varepsilon$ and $T$. 
\end{proof}

With Duhamel's formula, the Strichartz estimate stated in Theorem \ref{strichartz} and the estimates on $h_{\varepsilon}$ and $g_{\varepsilon}$ stated in Proposition \ref{estimateh} and Proposition \ref{estimateg} respectively, we obtain
\begin{equation}\nonumber
\begin{split}
    \norm{\chi_{\varepsilon}}_{L^{\infty}_{0\le t\le T}(L^2_x)} & = \norm{\int_0^t e^{-\mathrm{i}(t-\tau)H_x}(h_{\varepsilon} +g_{\varepsilon})(\tau,x) \mathrm{d} \tau }_{L^{\infty}_{0\le t\le T}(L^2_x)}  \quad \text{(Duhamel's formula)}  \\
    & \le C_{} \norm{h_{\varepsilon} +g_{\varepsilon}}_{L^{\frac{4}{3}}_{0\le t\le T}(L^1_x)} \quad \text{(Theorem \ref{strichartz})}\\
    & \le C_{}  \norm{h_{\varepsilon}}_{L^{\frac{4}{3}}_{0\le t\le T}(L^1_x)} + C_{}   \norm{ g_{\varepsilon}}_{L^{\frac{4}{3}}_{0\le t\le T}(L^1_x)} \\
    & \le C_{}  T^{\frac{3}{4}} \varepsilon^{\frac{1}{4}} + C_{}  \norm{\chi_{\varepsilon}}_{L^{\frac{4}{3}}_{0\le t\le T}(L^2_x)}. \quad \text{(Proposition \ref{estimateh} and \ref{estimateg})}
\end{split}
\end{equation}
It results in that, for $T\in [0,T_0]$,
\begin{equation} \label{beforeGronwall}
        \norm{\chi_{\varepsilon}(T,\cdot)}_{L^2}^{\frac{4}{3}} \le \norm{\chi_{\varepsilon}}_{L^{\infty}_{0\le t\le T}(L^2_x)}^{\frac{4}{3}} \le C_{}  T_0 \varepsilon^{\frac{1}{3}} + C_{}  \int_0^T\norm{\chi_{\varepsilon}(t,\cdot)}_{L^2}^{\frac{4}{3}} \mathrm{d }t.
\end{equation}
Applying Gr\"onwall's inequality in integral form~\cite[Appendix B.2]{Evans} to \eqref{beforeGronwall}, it follows that 
\begin{equation}\nonumber
    \norm{\chi_{\varepsilon}(T,\cdot)}_{L^2}^{\frac{4}{3}}  \le C T_0 \varepsilon^{\frac{1}{3}} \LR{1+C T e^{C T}}\le C_{T_0} \varepsilon^{\frac{1}{3}} 
\end{equation}
for a constant $C_{T_0}$ only depending on $T_0$, which implies Theorem \ref{difference}. \\

With direct calculations, we have
\begin{equation}\nonumber
    \norm{\phi_{\varepsilon}(t,\cdot) - e^{-\mathrm{i}\frac{t}{\varepsilon} } \varphi(t,\cdot)u_{1} }_{L^2(\mathbb{R}^2)} \le \norm{\LR{\phi_{\varepsilon} - \Pi_1 \phi_{\varepsilon}}(t,\cdot)}_{L^2(\mathbb{R}^2)} + \norm{\LR{\varphi_{\varepsilon} - \varphi}(t,\cdot)}_{L^2(\mathbb{R})},
\end{equation}

\begin{multline}\nonumber
    \norm{\psi(t,\cdot) - e^{-\mathrm{i}\frac{t}{\varepsilon} } \varphi(t,\cdot)u_{\varepsilon} }_{L^2(\mathbb{R}^2)} \le \norm{\phi_{\varepsilon}(t,\cdot) - e^{-\mathrm{i}\frac{t}{\varepsilon} } \varphi(t,\cdot)u_{1} }_{L^2(\mathbb{R}^2)} \\ +  \norm{ \LR{1 - e^{ -\mathrm{i} \beta S_{\varepsilon}[  |\phi_{\varepsilon}|^2](t,\cdot
    )} } \varphi(t,x)u_{1}(y) }_{L^2(\mathbb{R}^2)}.
\end{multline}
Applying Lemma \ref{diffphiestim}, Theorem \ref{difference} and the dominated convergence theorem, we conclude the proof of Theorem \ref{dynamicresult}.

\begin{appendix}

\section{Conservative vector fields and irrotational vector fields}

Let $\Omega$ be an open and simply connected subset of $\mathbb{R}^2$. Consider the vector field $\mathbf{A} : \Omega \to \mathbb{R}^2 $. In the following, we focus on the Sobolev space $ \mathbf{W}^{1,1}_{}(\Omega)$, and the equations and the derivatives hold in the sense of distribution (in the weak sense).

\begin{definition}[Conservative vector fields]
    We say that $\mathbf{A}$ is conservative if there exists a function $\Phi$ such that
    \begin{equation}\nonumber
        \mathbf{A} = \nabla \Phi.
    \end{equation}
\end{definition}

\begin{definition}[Irrotational vector fields]
    For $\mathbf{A} = (A_x, A_y)$, its Curl is defined as
    \begin{equation}\nonumber
        \operatorname{{Curl}} \mathbf{A} = \partial_x A_y - \partial_y A_x.
    \end{equation}

    We say that $\mathbf{A}$ is irrotational if its $\operatorname{Curl}$ is zero.
\end{definition}

\begin{theorem}[Equivalence]\label{0curl}
If $\mathbf{A} \in \mathbf{W}^{1,1}_{}(\Omega)$, then $\mathbf{A}$ is conservative if and only if it is irrotational.
\end{theorem}

\begin{proof}
{\noindent\bf{conservative $\Rightarrow$ irrotational.}}
\begin{equation*}
    \operatorname{{Curl}} \mathbf{A} = \partial_x \partial_y\Phi - \partial_y \partial_x\Phi = 0.
\end{equation*}

{\noindent\bf{irrotational $\Rightarrow$ conservative.}}

Fix $\mathbf{x}_0 \in \Omega$.
Define 
\begin{equation}\label{defPotential}
\Phi(\mathbf{x}) = \int_{\gamma_\mathbf{x}} \mathbf{A} \cdot \differential \mathbf{p}, \end{equation}
where $\gamma_{\mathbf{x}}$ is a nice enough curve contained in $\Omega$ starting from $\mathbf{x}_0$ and ending at $\mathbf{x}$. 
Since $\mathbf{A} \in \mathbf{W}^{1,1}_{}(\Omega)$, thanks to the trace theorem and the Sobolev embeddings, we have $\mathbf{A} \in \mathbf{L}^{1}_{}(\gamma_{\mathbf{x}})$, so the integral on the right-hand side of the definition \eqref{defPotential} is well-defined. 
Now we should show that the value of $\Phi(\mathbf{x})$ is independent of the choice of curve $\gamma_{\mathbf{x}}$. Let $\gamma_{\mathbf{x}}$ and $\gamma_{\mathbf{x}}'$ be two different nice enough curves contained in $\Omega$ starting from $\mathbf{x}_0$ and ending at $\mathbf{x}$, and let $D_{\gamma_{\mathbf{x}}, \gamma_{\mathbf{x}}'}$ be the surface enclosed by these two curves. With the help of Green's formula, we have the following:
\begin{equation*}
     \int_{\partial D_{\gamma_{\mathbf{x}}, \gamma_{\mathbf{x}}'}} \mathbf{A} \cdot \differential \mathbf{p} = \int_{D_{\gamma_{\mathbf{x}}, \gamma_{\mathbf{x}}'}} \operatorname{Curl} \mathbf{A} \, \differential x \differential y = 0,
\end{equation*}
which implies that 
\begin{equation*}
    \int_{\gamma_\mathbf{x}} \mathbf{A} \cdot \differential \mathbf{p} = \int_{\gamma_\mathbf{x}'} \mathbf{A} \cdot \differential \mathbf{p},
\end{equation*}
i.e. $\Phi$ is a well-defined function. The condition $\mathbf{A} \in \mathbf{W}^{1,1}_{}(\Omega)$ and the simple connectedness of $\Omega$ ensure the validity of Green's formula above. Then we are going to prove the gradient of $\Phi$ is equal to $\mathbf{A}$. With the direct calculations, for $(x,y)\in \Omega$ and small enough $h>0$, we have
\begin{equation*}
\begin{split}
    \frac{\Phi(x+h,y) - \Phi(x-h,y)}{2h} & = \frac{1}{2h} \int_{x-h}^{x+h} A_x(t,y) \differential t \\
    & = \frac{1}{2h}  \int_{\mathbb{R}} \mathbbm{1}_{[-h,h]}(x-t) (\mathbbm{1}_{\Omega}A_x)(t,y)\differential t \\
    & = \left(\frac{\mathbbm{1}_{[-h,h]}}{2h} * (\mathbbm{1}_{\Omega}A_x)(\cdot,y)\right)(x),
\end{split}
\end{equation*}
where $\mathbf{A} = (A_x, A_y)$ and $\mathbbm{1}$ is the indicator function. Clearly, $\left \{ \frac{\mathbbm{1}_{[-h,h]}}{2h}\right\}_{h>0}$ is an approximation to the identity, i.e. $\frac{\mathbbm{1}_{[-h,h]}}{2h}$ converges to the delta function $\delta_0$ as $h$ goes to 0. Since $ (\mathbbm{1}_{\Omega}A_x)(\cdot,y) \in L^1(\mathbb{R}) $, we have that $\frac{\mathbbm{1}_{[-h,h]}}{2h} * (\mathbbm{1}_{\Omega}A_x)(\cdot,y)$ converges strongly to $(\mathbbm{1}_{\Omega}A_x)(\cdot,y) $ in $L^1(\mathbb{R})$ as $h$ goes to 0. Therefore, 
\begin{multline}\nonumber
    \partial_x \Phi (x,y) = \lim_{h\to 0} \frac{\Phi(x+h,y) - \Phi(x-h,y)}{2h} = \lim_{h\to 0} \left(\frac{\mathbbm{1}_{[-h,h]}}{2h} * (\mathbbm{1}_{\Omega}A_x)(\cdot,y)\right)(x)  \\ = (\mathbbm{1}_{\Omega}A_x)(x,y) = A_x(x,y).
\end{multline}
Similarly, we have 
\begin{equation*}
    \partial_y \Phi = A_y.
\end{equation*}
Hence, 
\begin{equation*}
    \nabla \Phi = \mathbf{A}.
\end{equation*}
\end{proof}

\section{Continuity equation from Schr\"odinger equations}

In this section, we derive the continuity equation associated to 
\begin{equation}
    \mathrm{i} \partial_t \psi = \left( -\mathrm{i} \nabla_{\mathbf{x}} + \mathbf{A}(\psi) \right)^2 \psi  + g(\psi) \psi, \label{SchEq}
\end{equation}
where $\mathbf{A}(\psi) $ and $g(\psi) $ are both real-valued functions.\\

We denote 
\begin{equation}
    \rho = | \psi |^2. \nonumber
\end{equation}
From \eqref{SchEq} we have
\begin{gather}
    \partial_t \psi  = - \mathrm{i}  \left[ \left( -\mathrm{i} \nabla_{\mathbf{x}} + \mathbf{A}(\psi)  \right)^2 \psi + g(\psi)  \psi \right], \quad 
    \partial_t \overline{\psi} =  \mathrm{i}  \left[ \left( \mathrm{i} \nabla_{\mathbf{x}} + \mathbf{A}(\psi)  \right)^2 \overline{\psi} + g(\psi)  \overline{\psi} \right]. \nonumber
\end{gather}
Thus
\begin{align}
    \partial_t \rho &  = \overline{\psi} \partial_t \psi + \psi \partial_t \overline{\psi}  \nonumber \\
    & = - \mathrm{i} \overline{\psi} \left[ \left( -\mathrm{i} \nabla_{\mathbf{x}} + \mathbf{A}(\psi)  \right)^2 \psi + g(\psi)  \psi \right] + c. c. \nonumber \\
    & =   -\mathrm{i} \overline{\psi}  \left( -\mathrm{i} \nabla_{\mathbf{x}} + \mathbf{A}(\psi)  \right)^2 \psi  + c. c.   \nonumber \\
    & = - \mathrm{i} \overline{\psi}  \big[ - \mathrm{i} \nabla_{\mathbf{x}} \cdot \left( -\mathrm{i} \nabla_{\mathbf{x}} \psi + \mathbf{A}(\psi)  \psi \right)  +  ( -  \mathrm{i}\nabla_{\mathbf{x}} \psi) \cdot  \mathbf{A}(\psi)  \big] + c. c. \nonumber \\
    & =  - \overline{\psi} \nabla_{\mathbf{x}} \cdot \left( -\mathrm{i} \nabla_{\mathbf{x}} \psi + \mathbf{A}(\psi)  \psi \right) -  ( \nabla_{\mathbf{x}} \psi) \cdot  (\mathbf{A}(\psi) \overline{\psi} )  + c. c.\nonumber \\
    & =  - \overline{\psi} \nabla_{\mathbf{x}} \cdot \left( -\mathrm{i} \nabla_{\mathbf{x}} \psi + \mathbf{A}(\psi)  \psi \right) -  ( \nabla_{\mathbf{x}} \psi) \cdot  (\mathrm{i} \nabla_{\mathbf{x}} \overline{\psi} +  \mathbf{A}(\psi) \overline{\psi} )  \nonumber\\
    & \qquad  \qquad  \qquad  \qquad  \qquad  \qquad  - {\psi} \nabla_{\mathbf{x}} \cdot \left( \mathrm{i} \nabla_{\mathbf{x}} \overline{\psi} + \mathbf{A}(\psi)  \overline{\psi} \right) -  ( \nabla_{\mathbf{x}} \overline{\psi}) \cdot  (- \mathrm{i} \nabla_{\mathbf{x}} {\psi} + \mathbf{A}(\psi) {\psi} ) \nonumber \\
    & = - \nabla_{\mathbf{x}} \cdot \left[ \overline{\psi} \left( -\mathrm{i} \nabla_{\mathbf{x}} + \mathbf{A}(\psi)  \right) \psi  \right] - \nabla_{\mathbf{x}} \cdot \left[ {\psi} \left( \mathrm{i} \nabla_{\mathbf{x}} + \mathbf{A}(\psi)  \right) \overline{\psi} \right] \nonumber\\
    & = - \nabla_{\mathbf{x}} \cdot \left[ \overline{\psi} \left( -\mathrm{i} \nabla_{\mathbf{x}} + \mathbf{A}(\psi)  \right) \psi  + {\psi} \left( \mathrm{i} \nabla_{\mathbf{x}} + \mathbf{A}(\psi)  \right) \overline{\psi} \right]\nonumber\\
    & = - \nabla_{\mathbf{x}} \cdot \left[ \overline{\psi} \left( -\mathrm{i} \nabla_{\mathbf{x}} + \mathbf{A}(\psi)  \right) \psi  + {\psi}  \overline{ \left( -\mathrm{i} \nabla_{\mathbf{x}} + \mathbf{A}(\psi)  \right)\psi} \right] \nonumber\\
    & = - 2 \nabla_{\mathbf{x}} \cdot  \mathbf{J}_{\mathbf{A}(\psi) } (\psi), \nonumber
\end{align}
where $\mathbf{J}$ is as in \eqref{defJ}. Hence, the continuity equation associated to \eqref{SchEq} is 
\begin{equation}
    \partial_t \rho  + 2 \nabla_{\mathbf{x}} \cdot \mathbf{J}_{\mathbf{A}(\psi) } (\psi) = 0 . \label{cteq}
\end{equation}

\section{Equivalent norms for $\Sigma^s (\mathbb{R}^n) $}\label{equivnorm}

Recall Definition \ref{defSigma} and notice that $\Sigma^s (\mathbb{R}^n) $ is a Hilbert space with the equivalent norms
\begin{equation}\nonumber
    \norm{\Psi}_{\Sigma^s(\mathbb{R}^n)}^2 = \norm{\Psi}_{H^s(\mathbb{R}^n)}^2 + \norm{|\mathbf{x}|^s\Psi}_{L^2(\mathbb{R}^n)}^2 \simeq \norm{\Psi}_{L^2(\mathbb{R}^n)}^2 + \norm{(-\Delta + |\mathbf{x}|^2)^{\frac{s}{2}}\Psi}_{L^2(\mathbb{R}^n)}^2.
\end{equation}
The norm equivalence is obvious for $s=1$. For $s=2$, using integration by parts, we have
\begin{equation}\nonumber
    \norm{(-\Delta + |\mathbf{x}|^2)^{}\Psi }_{L^2(\mathbb{R}^n)}^2 = \norm{\Delta \Psi}_{L^2(\mathbb{R}^n)}^2  + \norm{|\mathbf{x}|^2\Psi}_{L^2(\mathbb{R}^n)}^2 - 2n\norm{\Psi}_{L^2(\mathbb{R}^n)}^2 + 2 \norm{\mathbf{x}\cdot\nabla\Psi}_{L^2(\mathbb{R}^n)}^2 ,
\end{equation}
which implies
\begin{equation}\nonumber
    \norm{\Psi}_{H^2(\mathbb{R}^n)}^2 + \norm{|\mathbf{x}|^2\Psi}_{L^2(\mathbb{R}^n)}^2 +\norm{\mathbf{x}\cdot\nabla\Psi}_{L^2(\mathbb{R}^n)}^2  \le  C \LR{\norm{\Psi}_{L^2(\mathbb{R}^n)}^2 + \norm{(-\Delta + |\mathbf{x}|^2)^{}\Psi}_{L^2(\mathbb{R}^n)}^2}
\end{equation}
for some constant $C$ independent of $\Psi$.
On the other hand, using triangle inequality, we find
\begin{equation}\nonumber
    \norm{(-\Delta + |\mathbf{x}|^2)^{}\Psi }_{L^2(\mathbb{R}^n)}^2 \le 2 \norm{\Delta \Psi}_{L^2(\mathbb{R}^n)}^2  + 2\norm{|\mathbf{x}|^2\Psi}_{L^2(\mathbb{R}^n)}^2,
\end{equation}
which shows 
\begin{equation}\nonumber
    \norm{\Psi}_{L^2(\mathbb{R}^n)}^2 + \norm{(-\Delta + |\mathbf{x}|^2)^{}\Psi}_{L^2(\mathbb{R}^n)}^2 \le C \LR{\norm{\Psi}_{H^2(\mathbb{R}^n)}^2 + \norm{|\mathbf{x}|^2\Psi}_{L^2(\mathbb{R}^n)}^2  }
\end{equation}
for some constant $C$ independent of $\Psi$.

For other values of $s$, the proof for this norm equivalence might be much more subtle to deal with; readers may refer to \cite[Section~2]{Ben-Cas-Meh} for example.
However, we only need the case $s = 1,2$ throughout the text.

\section{A Strichartz estimate}

\begin{definition}
    A pair $(p,q)$ is called {\it Strichartz-admissible} if 
    \begin{equation}\nonumber
        \frac{2}{p} +\frac{1}{q} = \frac{1}{2}, \quad q\in [2,\infty].
    \end{equation}
    For $a\in[1,\infty]$, we define $a'$ as the {\it conjugate exponnent} of $a$ by 
    \begin{equation}\nonumber
        \frac{1}{a} + \frac{1}{a'} = 1.
    \end{equation}
\end{definition}

\begin{theorem}[Strichartz estimates]\label{strichartz}
    Let $H_x$ be as in \eqref{H_xH_y}. For any Strichartz-admissible pairs $(p,q)$ and $(a,b)$ and for function $F :[0,T_0] \times \mathbb{R} \to \mathbb{C}$, there exists $C_{}>0$ only depending on $p,q,a,b,T_0$ such that 
    \begin{equation}\nonumber
        \norm{\int_0^t e^{-\mathrm{i}(t-\tau)H_x}F(\tau,x) \mathrm{d} \tau }_{L^{p}_{0\le t\le T}(L^q_x)} \le C_{} \norm{F}_{L^{a'}_{0\le t\le T}(L^{b'}_x)}, \quad \forall T\in [0,T_0].
    \end{equation}
\end{theorem}

Readers can refer to Section 2.3, Section 2.7 and Section 9.2 of \cite{Caz-03} for the proof of Theorem \ref{strichartz}. 
Section 2.3 provides the Strichartz estimates for $-\Delta$ instead of $H_x$, but, as stated in Section 2.7, with Lemma 9.2.4 in Section 9.2, the Strichartz estimates are also true for $H_x$. 
Notice that we only care about the case of dimension one in space. The case $q=\infty$ is not the endpoint case that causes trouble.
That is to say, a similar proof of Theorem 2.3.3 in Section 2.3 of \cite{Caz-03} is enough for Theorem \ref{strichartz}.
The endpoint case for higher-dimensional spaces can be found in \cite{KeeTao-98}.
    
\end{appendix}

\end{document}